\documentclass[10pt]{article}
\usepackage{amssymb, amsmath, amsfonts,latexsym, oldgerm,amscd,amsthm}
\usepackage{relsize}
\usepackage{url}

\textwidth13cm \usepackage[all]{xy} 
\newcommand{\tn}{\textnormal} \newcommand{\lra}{\longrightarrow}
 \newtheorem{de}{Definition}[section]

 \newtheorem{pr}[de]{Proposition}
\newtheorem{thm}[de]{Theorem} \newtheorem{lem}[de]{Lemma}
 \newtheorem{nt}[de]{Notation}
 \newtheorem{cor}[de]{Corollary}
 \newtheorem{rmk}[de]{Remark}

\def\E{{\rm E}}

\def\ESp{{\rm ESp}}
\def\Sp{{\rm Sp}}

\def\GL{{\rm GL}}
\def\G{{\rm G}}

\def\Um{{\rm Um}}

\def\End{{\rm End}}
\def\Aut{{\rm Aut}}
\def\Trans{{\rm Trans}}
\def\ETrans{{\rm ETrans}}
\def\sr{{\rm sr}}

\def\lra{\longrightarrow}

\newcommand{\gm}{\mathfrak{m}}

\begin{document}
\title{Equality of Linear and Symplectic Orbits}
\author{Pratyusha Chattopadhyay and
 Ravi A. Rao
  \\
{\small Stat-Math Unit, Indian Statistical Institute, 203 B.T. Road,
 Kolkata 700 108, India}
 {\small \&}\\ {\small Tata
    Institute of Fundamental Research, Dr. Homi Bhabha Road, Mumbai
    400 005, India}
} \date{}
\maketitle

\begin{center}{\it 2000 Mathematics Subject Classification: {13C10, 
13H05, 15A63, 19A13, 19B10, 19B14}}
\end{center}

\begin{center}{\it Key words: unimodular rows, transvections, elementary symplectic group}
\end{center} {\small ~~~~Abstract: It is shown that the set of orbits
  of the action of the elementary symplectic transvection group on all
  unimodular elements of a symplectic module over a commutative ring
  in which 2 is invertible is identical with the set of orbits of the
  action of the corresponding elementary transvection group. This
  result is used to get improved injective stability estimates for
  $K_1$ of the symplectic transvection group over a non-singular
  affine algebras.}

\section{\large Introduction}

In this paper we discuss two related questions about (linear and
symplectic) elementary transvections of a projective $R$-module
(resp. symplectic $R$-module) of type $R \oplus P$
(resp. $(\mathbb{H}(R) \oplus (P, \langle , \rangle))$).

The first one is about comparing the groups generated by the two
different types of elementary transvections; and showing that they are
the same. (This fact seems to have escaped notice earlier; and experts
have told us that it is interesting, and opens up the study done on
these transvections.)

The second one is to show that the linear and symplectic elementary
transvection orbits of a unimodular element in a symplectic module
coincide. (This generalizes the result in \cite{cr} where it was shown
in the free case.)

We now describe the two problems a bit more in detail. 

H. Bass introduced two types of linear transvections of a projective
module $R \oplus P$ in \cite{bass}. He also introduced two types of
symplectic transvecitons of a symplectic module $\mathbb{H}(R) \oplus
(P, \langle, \rangle)$. (These are recalled in \S 4 and \S 5).

Since elementary automorphisms are homotopic to the identity, we are
able to invoke Quillen-Suslin theory (see \cite{Q}, \cite{18}) to show
that
\begin{itemize}

\item{the groups generated by the two types of elementary linear
  transvections are the same as the elementary linear group in the
  free case (see Lemma \ref{linear,free,rel}).}

\item{the group generated by the two types of elementary symplectic
  transvections w.r.t. the standard alternating form are the same as
  the elementary symplectic group in the free case (see Lemma
  \ref{free,psi}).}

\end{itemize}

The above generalizes the special case of these results in (\cite{BR},
Theorem 2). 

The title of this paper alludes to the comparison of the elementary
linear and elementary symplectic orbits of a unimodular element
$(a,b,p)$ in a symplectic module $(\mathbb{H}(R) \oplus P)$. In case
$P$ is free of rank $\geq 4$, it is established in (\cite{cr}, Theorem
4.2, Theorem 5.6) that these two orbits coincide. In the Appendix the
missing case when $P$ is free of rank $2$ is proved by a similar, but
slightly more involved argument. (This means one has to essentially
prove Lemma 2.9 and Lemma 3.1 in \cite{cr}.)

Since elementary automorphisms are homotopic to the identity, we show
how the Quillen-Suslin machinery enables one to extend the results of
\cite{cr} to show that the two orbits are equal in the general case
when $P$ is a finitely generated projective module. (The transition is
by no means automatic!)

\section{\large Preliminaries}

\medskip
A row ${ v} = (v_{1}, \ldots, v_{n}) \in R^{n}$ is said to be {\it
  unimodular} if there are elements $w_{1}, \ldots, w_{n}$ in $R$ such
that $v_{1}w_{1} + \cdots + v_{n}w_{n} = 1$. Um$_{n}(R)$ will denote
the set of all unimodular rows ${ v} \in R^{n}$. Let $I$ be an ideal
in $R$. We denote by ${\rm Um}_n(R, I)$ the set of all unimodular rows
of length $n$ which are congruent to $e_1 = (1, 0, \ldots, 0)$ modulo
$I$. (If $I = R$, then ${\rm Um}_n(R, I)$ is ${\rm Um}_n(R)$).

\begin{de} {\rm Let $P$ be a finitely generated projective $R$-module. 
An element $p \in P$ is said to be {\it unimodular} if there exists a 
$R$-linear map $\phi: P \to R$ such that $\phi(p)=1$. The collection 
of unimodular elements of $P$ is denoted by $\Um(P)$. 

Let $P$ be of the form $R \oplus Q$ and have an element of the form
$(1,0)$ which correspond to the unimodular element. An element $(a,q)
\in P$ is said to be {\it relative unimodular} w.r.t. an ideal $I$ of
$R$ if $(a,q)$ is unimodular and $(a,q)$ is congruent to $(1,0)$
modulo $IP$.  The collection of all relative unimodular elements
w.r.t. an ideal $I$ is denoted by ${\rm Um}(P,IP)$.  }
\end{de}

Let us recall that if $M$ is a finitely presented $R$-module and $S$
is a multiplicative set of $R$, then $S^{-1} {\rm Hom}_R(M,R) \cong
{\rm Hom}_{R_S}(M_S, R_S)$. Also recall that if $f=(f_1, \ldots,
f_n)\in R^n := M$, then $\Theta_M(f)=\{ \phi(f): \phi \in {\rm
  Hom}(M,R) \}= \sum_{i=1}^n Rf_i$. Therefore, if $P$ is a finitely
generated projective $R$-module of rank $n$, $\gm$ is a maximal ideal
of $R$ and $v\in \Um(P)$, then $v_\gm \in \Um_n(R_\gm)$. Similarly if
$v \in \Um(P,IP)$ then $v_\gm \in \Um_n(R_\gm, I_\gm)$.

The group GL$_{n }(R)$ of invertible matrices acts on $R^{n }$ in a
natural way: $ { v} \longrightarrow { v} \sigma$, if ${ v} \in R^{n}$,
$\sigma \in$ GL$_{n}(R)$. This map preserves Um$_{n}(R)$, so
GL$_{n}(R)$ acts on Um$_{n}(R)$. Note that any subgroup G of GL$_n(R)$
also acts on Um$_n(R)$. Let $v, w \in {\rm Um}_n(R)$, we denote $v
\displaystyle\sim_{{\rm G}} w$ or $v \in w$G if there is a $g \in$ G
such that $v= wg$.

Let E$_{n}(R)$ denote the subgroup of SL$_{n}(R)$ consisting of all
{\it elementary} matrices, i.e. those matrices which are a finite
product of the {\it elementary generators} E$_{ij}(\lambda) = I_{n} +
e_{ij}(\lambda)$, $1 \leq i \neq j \leq n$, $\lambda \in R$, where
$e_{ij}(\lambda) \in$ M$_{n}(R)$ has an entry $\lambda$ in its $(i,
j)$-th position and zero elsewhere.

In the sequel, if $\alpha$ denotes an $m \times n$ matrix, then we let
$\alpha^t$ denote its {\it transpose} matrix. This is of course an
$n\times m$ matrix. However, we will mostly be working with square
matrices, or rows and columns.

\begin{de} {\bf The Relative Groups E$_n(I)$, E$_n(R,I)$:} Let $I$ be
  an ideal of $R$. The relative elementary group {\rm E}$_n(I)$ is the
  subgroup of {\rm E}$_n(R)$ generated as a group by the elements {\rm
    E}$_{ij}(x)$, $x \in I$, $1 \leq i \neq j \leq n$.

  The relative elementary group ${\rm E}_n(R, I)$ is the normal
  closure of {\rm E}$_n(I)$ in {\rm E}$_n(R)$.

$($Equivalently, ${\rm E}_n(R, I)$ is generated as a group by ${\rm
    E}_{ij}(a) {\rm E}_{ji}(x)$E$_{ij}(-a)$,\tn{ with} $a \in R$, $x
  \in I$, $i \neq j$, provided $n \geq 3$ {\rm (see \cite{SV}, Section
    2)}$)$.
\end{de}

\begin{de} 
{\rm {\bf E$_n^1(R, I)$} is the subgroup of $\E_n(R)$ generated by the
  elements of the form $E_{1i}(a)$ and $E_{j1}(x)$, where $a \in R, x
  \in I$, and $2 \le i, j \le n$.  }
\end{de}

\begin{rmk} \label{unimodularOverLocalRing}
It is easy to check that if $v \in \Um_n(R,I)$, where $(R, \gm)$ is a
local ring and $I$ be an ideal of $R$, then $v = e_1 \beta$, for some
$\beta \in \E_n(R,I)$.
\end{rmk}

\begin{de} {\bf Symplectic Group Sp$_{2n}(R)$:} The symplectic group ${\rm
  Sp}_{2n}(R)=\{\alpha \in {\rm GL}_{2n}(R)\,\,|\,\, \alpha^t \psi_n
  \alpha = \psi_n \}$, where $\psi_n = \underset{i=1}{\overset{n}\sum}
  e_{2i-1,2i}- \underset{i=1}{\overset{n}\sum} e_{2i,2i-1}$, the
  standard symplectic form.
\end{de}

\tn{Let $\sigma$ denote the permutation of the natural numbers given
  by $\sigma(2i)=2i-1$ and $\sigma(2i-1)=2i$}.

\begin{de} \label{2.4} {\bf Elementary Symplectic Group
    ESp$_{2n}(R)$:} We define for $z \in R$, $1\le i\ne j\le 2n$,

\[ se_{ij}(z) =
\begin{cases}
  1_{2n}+ze_{ij}  & \mbox{{\rm if}  $i=\sigma(j)$,} \\
  1_{2n}+ze_{ij}-(-1)^{i+j}z e_{\sigma (j) \sigma (i)} & \mbox{{\rm
      if} $i\ne \sigma(j)$.}
\end{cases}
\]

\tn{(It is easy to check that all these elements belong to
  {\rm Sp}$_{2n}(R)$)}.

We call them elementary symplectic matrices over $R$ and the
subgroup of {\rm Sp}$_{2n}(R)$ generated by them is called the
elementary symplectic group {\rm ESp}$_{2n}(R)$.
\end{de}

\begin{de} {\bf The Relative Group ESp$_{2n}(I)$, ESp$_{2n}(R,I)$:}
  Let $I$ be an ideal of $R$. The relative elementary group ${\rm
  ESp}_{2n}(I)$ is the subgroup of ${\rm ESp}_{2n}(R)$ generated as a
  group by the elements $se_{ij}(x)$, $x \in I$ and $1 \le i \ne j \le
  2n$.

The relative elementary group ${\rm ESp}_{2n}(R,I)$ is the normal
    closure of ${\rm ESp}_{2n}(I)$ in ${\rm ESp}_{2n}(R)$.
\end{de}

\begin{de}
  The group $\ESp_{2n}^1(R,I)$ is the subgroup of $\ESp_{2n}(R)$
  generated by the elements of the form $se_{1 i}(a)$ and $se_{j
    1}(x)$, where $a \in R$, $x \in I$ and $2 \le i, j \le 2n$.
\end{de}

\begin{nt} {\rm We fix some notations. ${\rm G}(n,R)$ will denote 
either the linear group ${\rm GL}_n(R)$, or the symplectic group ${\rm
  Sp}_{2m}(R)$, {\bf for} {\bf $n=2m$}. Also, ${\rm E}(n,R)$ will
denote either of the elementary subgroups ${\rm E}_n(R)$ or ${\rm
  ESp}_{2m}(R)$. The standard elementary generators of ${\rm E}(n,R)$
are denoted by $ge_{ij}(a)$, $a\in R$. Let $I$ be an ideal in $R$. Let
$\G(n, R, I)$ denote the relative linear groups $\GL_n(R,I)$, or the
relative symplectic group $\Sp_{2m}(R,I)$. Note that ${\rm E}(n,I)$ is
a subgroup of ${\rm E}(n,R)$ generated as a group by the elements
$ge_{ij}(x)$, where $x \in I$, and $1 \le i \ne j \le n$.  Here, ${\rm
  E}(n,R,I)$ denotes the corresponding relative elementary subgroups
${\rm E}_n(R,I)$, ${\rm ESp}_{2m}(R,I)$, respectively.  These are the
normal closures of the subgroups ${\rm E}(n,I)$ in ${\rm E}(n,R)$,
which are also known to be generated by the elements
$ge_{ij}(a)ge_{ji}(x)ge_{ij}(-a)$, $a\in R,x\in I$, and $1 \le i \ne j
\le n$ (see \cite{V3}, Lemma 8). Also, ${\rm E}^1(n,R,I)$ is a
subgroup of $\E(n,R)$, generated by the elements $ge_{1i}(a)$, where
$a \in R$ and $ge_{j1}(x)$, where $x \in I$, $2 \le i, j \le n$ in the
linear and symplectic case.}
\end{nt}

\begin{lem} \label{vanderk1} 
Let $n \ge 3$ in the linear case and $n \ge 4$ in the symplectic
case. Let $R$ be a commutative ring with $R=2R$, and let $I$ be an
ideal of $R$. Then the following sequence is exact

\[ 1 \lra \E(n,R,I) \lra \E^1(n,R,I) \lra \E^1(n,R/I,0) \lra 1. \]

Thus $\E(n,R,I)$ equals $\E^1(n,R,I) \cap \G(n,R,I)$.
\end{lem}

Proof: For the proof in the linear case reader can see (\cite{vdK1},
Lemma 2.2). The proof in the symplectic case is similar when $n \geq
6$. The proof when $n = 4$ is discussed in the Appendix.
\hfill{$\square$}

\section{{\large Local Global Principle for Relative Elementary Group}}

In this section we reprove the Local-Global Principle for an extended
ideal, done in \cite{acr}, by a different approach. The idea in
\cite{acr} is to use the {\bf special forms} (as described by Suslin
in \cite{18} for the linear group, Kopeiko in \cite{29} for the
symplectic group, and Suslin-Kopeiko in \cite{SK} for the orthogonal
groups). The idea here is to use a slightly larger group E$^1(n, R,
I)$ than the relative groups E$(n, R, I)$.  This group was introduced
by W. van der Kallen in \cite{vdK1}. Note that the absolute case (when
$I = R$) was first done in \cite{bbr} again by a different approach.
We proceed to describe this argument now.

\begin{lem} \label{ness4-E1} Let $R$ be a commutative ring with
  $R=2R$, and let $I$ be an ideal of $R$. Let $n \ge 3$ in the linear
  case and $n \ge 4$ in the symplectic case. Let $\varepsilon =
  \varepsilon_1 \ldots \varepsilon_r$ be an element of $\E^1(n,R,I)$,
  where each $\varepsilon_k$ is an elementary generator.  If
  $ge_{ij}(Xf(X))$ is an elementary generator of $\E^1(n,R[X],I[X])$,
  then
\begin{eqnarray*}
\varepsilon ~ ge_{ij}(Y^{4^r}Xf(Y^{4^r}X)) ~ \varepsilon^{-1} 
&=& \prod_{t=1}^s ge_{i_t j_t}(Y h_t(X,Y)),
\end{eqnarray*}
where either $i_t=1$ or $j_t=1$ and $h_t(X,Y) \in R[X,Y]$, when
$i_t=1$; $h_t(X,Y) \in I[X,Y]$ when $j_t=1$.
\end{lem}

Proof: Proof in the linear case is straight forward. We establish the
result here in the symplectic case when $n \geq 6$. In the Appendix
the case when $n = 4$ is also covered.

We prove the result using induction on $r$, where $\varepsilon$ is
product of $r$ many elementary generators. Let $r=1$ and $\varepsilon
= se_{pq}(a)$. Note that $a \in R$ when $p=1$, and $a \in I$ when
$q=1$. Also, note that we use $*$ to represent elements of the ideal
$I$. Given that $se_{ij}(Xf(X))$ is an elementary generator of
$\ESp_{2n}^1(n,R[X],I[X])$. First we assume $i=1$ and $j \ne 2$, hence
$f(X) \in R[X]$.

{\it Case $($1$)$:} Let $(p,q)=(1,j)$. In this case
\begin{eqnarray*}
se_{1j}(a) ~ se_{1j}(Y^4X f(Y^4X)) ~ se_{1j}(-a) &=& se_{1j}(Y^4X f(Y^4X)).
\end{eqnarray*}

{\it Case $($2$)$:} Let $(p,q) = (1, \sigma(j))$. In this case
\begin{eqnarray*}
&& se_{1 \sigma(j)}(a) ~ se_{1j}(Y^4X f(Y^4X)) ~ se_{1 \sigma(j)}(-a) \\ 
&=& se_{12}(-a Y^4X f(Y^4X)) ~ se_{1j}(Y^4X f(Y^4X)).
\end{eqnarray*}

{\it Case $($3$)$:} Let $(p,q)=(1,k), k \ne j, \sigma(j)$. In this case
\begin{eqnarray*}
se_{1k}(a) ~ se_{1j}(Y^4X f(Y^4X)) ~ se_{1k}(-a) &=& se_{1j}(Y^4X f(Y^4X)).
\end{eqnarray*}

{\it Case $($4$)$:} Let $(p,q)=(2,1)$. In this case
\begin{eqnarray*}
&& se_{21}(a) ~ se_{1j}(Y^4X f(Y^4X)) ~ se_{21}(-a) \\
&=& se_{2j}(*Y^4X f(Y^4X)) ~ se_{\sigma(j) j}(*Y^4X f(Y^4X)) ~ se_{1j}(Y^4X f(Y^4X)) \\
&=& se_{\sigma(j) 1}(*Y^4X f(Y^4X)) ~ [se_{\sigma(j) 1}(*Y^2), ~ se_{1j}(Y^2X f(Y^4X))] \\
&& se_{1j}(Y^4X f(Y^4X)) .
\end{eqnarray*}

{\it Case $($5$)$:} Let $(p,q)=(k,1), k \ne 2, j, \sigma(j)$. In this case
\begin{eqnarray*}
&& se_{k1}(a) ~ se_{1j}(Y^4X f(Y^4X)) ~ se_{k1}(-a) \\
&=& se_{kj}(*Y^4Xf(Y^4X)) ~ se_{1j}(Y^4X f(Y^4X))\\ 
&=& [se_{k1}(*Y^2), ~ se_{1j}(Y^2Xf(Y^4X))] ~ se_{1j}(Y^4X f(Y^4X)).
\end{eqnarray*}

{\it Case $($6$)$:} Let $(p,q)=(\sigma(j),1)$. In this case
\begin{eqnarray*}
&& se_{\sigma(j) 1}(a) ~ se_{1j}(Y^4X f(Y^4X)) ~ se_{\sigma(j) 1}(-a) \\
&=& se_{\sigma(j) j}(*Y^4X f(Y^4X)) ~ se_{1j}(Y^4X f(Y^4X)) \\
&=& [se_{\sigma(j) 1}(*Y^2), ~ se_{1j}(Y^2X f(Y^4X))] ~ se_{1j}(Y^4X f(Y^4X)) .
\end{eqnarray*}

{\it Case $($7$)$:} Let $(p,q)=(j,1)$. Let us choose an integer $d \ne
1, 2, j, \sigma(j)$. In this case
\begin{eqnarray*}
&& se_{j 1}(a) ~ se_{1 j}(Y^4 Xf(Y^4 X)) ~ se_{j 1}(-a) \\ 
&=&  se_{j 1}(a) ~ [se_{1 d}(Y^2 Xf(Y^4 X)), ~ se_{d j}(Y^2)] ~ se_{j 1}(-a) \\
&=&  [se_{j d}(* Y^2 X f(Y^4 X)) ~ se_{1 d}(Y^2 X f(Y^4 X)), ~ se_{d 1}(* Y^2) 
 se_{d j}(Y^2)] \\
&=& se_{j d}(* Y^2 X f(Y^4 X)) ~ se_{1 d}(Y^2 X f(Y^4 X)) ~ se_{d 1}(* Y^2) ~ se_{d j}(Y^2) \\ 
&&  se_{1 d}(-Y^2 X f(Y^4 X)) ~ se_{j d}(- * Y^2 X f(Y^4 X)) ~ se_{d j}(- Y^2) 
 se_{d 1}(- * Y^2) \\
&=& se_{j d}(* Y^2 X f(Y^4 X)) ~ se_{1 d}(Y^2 X f(Y^4 X)) ~ se_{d 1}(* Y^2) ~ se_{d j}(Y^2) \\ 
&& se_{1 d}(-Y^2 X f(Y^4 X)) ~ [se_{j 1}(- * Y), se_{1 d}(Y X f(Y^4 X))] ~ se_{d j}(- Y^2)  
 se_{d 1}(- * Y^2) \\
&=& se_{j d}(* Y^2 X f(Y^4 X)) ~ se_{1 d}(Y^2 X f(Y^4 X)) ~ se_{d 1}(* Y^2) ~ se_{d j}(Y^2) \\ 
&& se_{1 d}(-Y^2 X f(Y^4 X)) ~ se_{d j}(-Y^2) ~ se_{d j}(Y^2) \\ 
&&  [se_{j 1}(-* Y), ~ se_{1 d}(Y X f(Y^4 X))] ~ se_{d j}(- Y^2) ~ se_{d 1}(- * Y^2) \\
&=&  [se_{j 1}(* Y), ~ se_{1 d}(Y X f(Y^4 X))] ~ se_{1 d}(Y^2 X f(Y^4 X)) ~ se_{d 1}(* Y^2) \\ 
&&  se_{1 d}(-Y^2 X f(Y^4 X)) ~ se _{1 j}(Y^4 X f(Y^4 X)) ~ [se_{d 1}(- * Y^3) se_{j 1}(- * Y), \\
&&  se_{1 j}(-Y^3 X f(Y^4 X)) ~ se_{1 d}(Y X f(Y^4 X))] ~ se_{d 1}(* Y^2).
\end{eqnarray*}

Now we consider the case when $(i, j)=(1, 2)$.

{\it Case $($8$)$:} Let $(p,q)=(1, 2)$. In this case
\begin{eqnarray*}
se_{12}(a) ~ se_{12}(Y^4X f(Y^4X)) ~ se_{12}(-a) &=& se_{12}(Y^4X f(Y^4X).
\end{eqnarray*}

{\it Case $($9$)$:} Let $(p,q)=(1, k)$, $k \ne 2$. In this case
\begin{eqnarray*}
se_{1k}(a) ~ se_{12}(Y^4X f(Y^4X)) ~ se_{1k}(-a) &=& se_{12}(Y^4X f(Y^4X).
\end{eqnarray*}

{\it Case $($10$)$:} Let $(p,q)=(k, 1)$, $k \ne 2$. In this case
\begin{eqnarray*}
&& se_{k1}(a) ~ se_{12}(Y^4X f(Y^4X)) ~ se_{k1}(-a) \\
&=& se_{1 \sigma(k)}(* Y^4X f(Y^4X)) ~ se_{k \sigma(k)}(* Y^4X f(Y^4X)) ~ se_{12}(Y^4X f(Y^4X))\\
&=& se_{1 \sigma(k)}(* Y^4X f(Y^4X)) ~ [se_{k1}(* Y^2), ~ se_{1 \sigma(k)}(Y^2X f(Y^4X))] \\
&& se_{12}(Y^4X f(Y^4X)).
\end{eqnarray*}

{\it Case $($11$)$:} Let $(p,q)=(2, 1)$. Let us choose an integer $d
\ne 1, 2$. In this case
\begin{eqnarray*}
&& se_{21}(a) ~ se_{12}(Y^4X f(Y^4X)) ~ se_{21}(-a) \\
&=& se_{21}(a) ~ [se_{1 d}(Y^2/2), ~ se_{d 2}(Y^2X f(Y^4X)) ] ~ se_{21}(-a) \\
&=& [se_{2 d}(* Y^2) ~ se_{\sigma(d) d}(* Y^4) ~ se_{1 d}(Y^2/2), se_{2 \sigma(d)}(* Y^2X f(Y^4X)) \\
&& se_{d \sigma(d)}(* Y^4X^2 f^2(Y^4X)) ~ se_{d 2}(Y^2X f(Y^4X)) ] \\
&=& [se_{\sigma(d) 1}(* Y^2) ~ [se_{\sigma(d) 1}(* Y^2), ~ se_{1 d}(Y^2)] ~ se_{1 d}(Y^2/2), se_{d 1}(* Y^2X f(Y^4X)) \\
&& [se_{d 1}(* Y^2), ~ se_{1 \sigma(d)}(Y^2X^2 f^2(Y^4X))] ~ se_{1 \sigma(d)}(Y^2X f(Y^4X)) ] 
\end{eqnarray*}

Hence the result is true when $i=1$ and $\varepsilon$ is an elementary
generator. Carrying out similar calculation one can show the result is
true when $j=1$ and $\varepsilon$ is an elementary generator. Let us
assume that the result is true when $\varepsilon$ is product of $r-1$
many elementary generators, i.e, $\varepsilon_2 \ldots \varepsilon_r ~
se_{ij}(Y^{4^{r-1}}X f(Y^{4^{r-1}}(X)) ~ \varepsilon_r^{-1} \ldots
\varepsilon_2^{-1} = \prod_{t=1}^k se_{p_t q_t}(Y g_t(X,Y))$, where
either $p_t=1$ or $q_t=1$. Note that $g_t(X,Y) \in R[X,Y]$ when
$p_t=1$ and $g_t(X,Y) \in I[X,Y]$ when $q_t=1$.

We now establish the result when $\varepsilon$ is product of $r$ many
elementary generators. We have
\begin{eqnarray*}
&& \varepsilon ~ se_{ij}(Y^{4^r}X f(Y^{4^r}X)) ~ \varepsilon^{-1}  \\
& = & \varepsilon_1 \varepsilon_2 \ldots \varepsilon_r ~ se_{ij}(Y^{4^r}X f(Y^{4^r}X)) 
~ \varepsilon_r^{-1} \ldots \varepsilon_2^{-1} \varepsilon_1^{-1} \\
& = & \varepsilon_1 ~ \big( \prod_{t=1}^k se_{p_t q_t}(Y^4 g'_t(X,Y)) \big) ~ \varepsilon_1^{-1}  \\
& = & \prod_{t=1}^k \varepsilon_1 ~ se_{p_t q_t}(Y^4 g'_t(X,Y)) ~ \varepsilon_1^{-1} \\
& = & \prod_{t=1}^s se_{i_t j_t}(Y h_t(X,Y)). 
\end{eqnarray*}

To get the last equality one needs to repeat the calculation which was
done for a single elementary generator. Note that at the last line
either $i_t=1$ or $j_t=1$. Also, note that $h_t(X,Y) \in R[X,Y]$, when
$i_t=1$ and $h_t(X,Y) \in I[X,Y]$, when $j_t=1$ 
\hfill{$\square$}

\begin{nt}
{\rm Let $M$ be a finitely presented $R$-module and $a$ be a
  non-nilpotent element of $R$. Let $R_a$ denote the ring $R$
  localized at the multiplicative set $\{a^i : i \ge 0 \}$ and $M_a$
  denote the $R_a$-module $M$ localized at $\{a^i : i \ge 0 \}$. Let
  $\alpha(X)$ be an element of $\End(M[X])$. The localization map $i:
  M \to M_a$ induces a map $i^*: \End(M[X]) \to \End(M[X]_a)
  = \End(M_a[X])$. We shall denote $i^*(\alpha(X))$ by $\alpha(X)_a$
  in the sequel.

}
\end{nt}

\begin{lem} \label{equal-auto}
  Let $M$ be a finitely presented $R$-module and $I$ be an ideal of
  $R$. Let $\alpha(X), \beta(X) \in
  {\End}(M[X],IM[X])=ker({\End}(M[X]) \lra {\End}(M[X]/IM[X]))$, with
  $\alpha(0)=\beta(0).$ Let $a$ be a non-nilpotent element in
  $R$ such that $\alpha(X)_a = \beta(X)_a$ in
  ${\End}(M_a[X],IM_a[X])$. Then $\alpha(a^N X)=\beta(a^N X)$ in
  ${\End}(M[X],IM[X])$, for $N \gg 0$.
\end{lem} 

Proof: Using $\alpha(0)-\beta(0)=0$, we get $\alpha(X)-\beta(X)=X
\gamma(X)$, for some $\gamma(X)$ in ${\rm End}(M[X],IM[X])$. Also
$\alpha(X)_a - \beta(X)_a = 0$ in ${\rm End}(M_a[X],IM_a[X])$, i.e,
$\big(\alpha(X) - \beta(X) \big)_a = 0$, i.e, $\big( X \gamma(X)
\big)_a = 0$. Hence $a^N \big( X \gamma(X) \big) = 0$, in ${\rm
  End}(M[X],IM[X])$, for some $N \gg 0$. Therefore
$\alpha(a^NX)-\beta(a^NX) = a^NX \gamma(a^NX) = 0$, in ${\rm
  End}(M[X],IM[X])$, for $N \gg 0$.  \hfill{$\square$}

\begin{lem} \label{E1-dil-strong} Let $R$ be a commutative ring with
  $R=2R$, and let $I$ be an ideal of $R$. Let $n \ge 3$ in the linear
  case and $n \ge 4$ in the symplectic case. Let $a$ be a
  non-nilpotent element in $R$ and $\alpha(X)$ be in
  $\E^1(n,R_a[X],I_a[X])$, with $\alpha(0)=Id$. Then there exists
  $\alpha^*(X) \in \E^1(n,R[X],I[X])$, with $\alpha^*(0) = Id.$, such
  that $\alpha^*(X)$ localises to $\alpha(bX)$, for $b \in (a^d)$, $d
  \gg 0$.
\end{lem}

Proof: Let $\alpha(X) = \prod_{k=1}^r ge_{i_k j_k}(f_k(X))$, where
$f_k(X)=f_k(0)+Xg_k(X)$. Therefore, $\alpha(X) = \prod_{k=1}^r
\gamma_k ~ ge_{i_k j_k}(Xg_k(X)) ~ \gamma_k^{-1}$, where $\gamma_l =
\prod_{k=1}^l ge_{i_k j_k}(f_k(0))$ and $\gamma_l$ is in
$\E^1(n,R_a,I_a)$. We can write $\alpha(Y^{4^r}X)=\prod_{k=1}^r \Big(
\prod_{t=1}^s ge_{i_t j_t}(Y h_t(X,Y))/a^{d'}) \Big)$, where either
$i_t=1$ or $j_t=1$ (see Lemma \ref{ness4-E1}). Note that $h_t(X,Y) \in
R[X,Y]$, when $i_t=1$ and $h_t(X,Y) \in I[X,Y]$, when $j_t=1$, and
$d'$ is a natural number. Let us choose $d=d'$ and define
$\alpha^*(X,Y)$ to be $\prod_{k=1}^r \Big( \prod_{t=1}^s ge_{i_t
  j_t}(Y h_t(X,a^dY)) \Big)$.

Clearly $\alpha^*(X,Y) \in \E^1(n,R[X,Y],I[X,Y])$ and $\alpha((a^d
Y)^{4^r} X) = \alpha^*(X,Y)$. Substituting $Y=1$, we get $\alpha(b
X) = \alpha^*(X)$, for $b \in (a^d)$, $d \gg 0$. Note that
$\alpha^*(X) \in \E^1(n,R[X],I[X])$, with $\alpha^*(0) =Id$.
\hfill{$\square$}

\begin{thm} \label{rel-dil-strong} Let $R$ be a commutative ring with
  $R=2R$, and let $I$ be an ideal of $R$. Let $n \ge 3$ in the linear
  case and $n \ge 4$ in the symplectic case. Let $a$ be a
  non-nilpotent element in $R$ and $\alpha(X)$ be in
  $\E(n,R_a[X],I_a[X])$, with $\alpha(0)=Id$. Then there exists
  $\alpha^*(X) \in \E(n,R[X],I[X])$, with $\alpha^*(0) = Id.$, such that
  $\alpha^*(X)$ localises to $\alpha(bX)$, for $b \in (a^d)$, $d \gg
  0$.
\end{thm}

Proof: Follows from the previous lemma and Lemma \ref{vanderk1},
which says that $\E(n,R,I)$ is $\E^1(n,R,I) \cap \G(n,R,I)$.
\hfill{$\square$}

\section{\large{Transvection Groups}}

Following H. Bass one  defines a transvection of a finitely generated
$R$-module as follows: 

\begin{de}
{\rm Let $M$ be a finitely generated $R$-module. Let $q \in M$ and
  $\pi \in M^*= {\rm Hom}(M,R)$, with $\pi(q) = 0$. Let $\pi_q(p) :=
  \pi(p)q$. An automorphism of the form $1+\pi_q$ is called a {\bf
    transvection} of $M$, if either $q \in \Um(M)$ or $\pi \in
  \Um(M^*)$. Collection of transvections of $M$ is denoted by
  $\Trans(M)$. This forms a subgroup of $\Aut(M)$.}
\end{de}

\begin{de}
{\rm Let $M$ be a finitely generated $R$ module. The automorphisms of
  $N= (R \oplus M)$ of the form
\begin{eqnarray*}
(a, p) & \mapsto & (a, p+ax),
\end{eqnarray*}
or of the form
\begin{eqnarray*}
(a, p) & \mapsto & (a+ \tau(p), p),
\end{eqnarray*}
where $x \in M$ and $\tau \in M^*$ are called {\bf elementary
  transvections} of $N$. Let us denote the first automorphism by $E_x$
and the second one by $E^*_\tau$.  It can be verified that these are
transvections of $N$.  Indeed, let us consider $\pi(t,y)=t$, $q=(0,x)$ to
get $E_x$, and consider $\pi(a,p)= \tau(p)$, where $\tau \in
M^*$, $q=(1,0)$ to get $E^*_\tau$. The subgroup of $\Trans(N)$
generated by elementary transvections is denoted by $\ETrans(N)$.}
\end{de}

\begin{de} {\rm Let $I$ be an ideal of $R$. The group of {\bf relative
      transvections} w.r.t. an ideal $I$ is generated by the
    transvections of the form $1+\pi_q$, where either $q \in IM$
    or $\pi \in IM^*$. The group generated by relative
    transvections is denoted by $\Trans(M,IM)$.}
\end{de}

\begin{de} {\rm Let $I$ be an ideal of $R$. The elementary
    transvections of $N=(R \oplus M)$ of the form $E_x, E^*_\tau$,
    where $x \in IM$ and $\tau \in IM^*$ are called {\bf relative
      elementary transvections} w.r.t. an ideal $I$, and the
    group generated by them is denoted by $\ETrans(IN)$. The normal
    closure of $\ETrans(IN)$ in $\ETrans(N)$ is denoted by
    $\ETrans(N,IN)$.}
\end{de}

\begin{lem} \label{linear,free,rel} Let $I$ be an ideal of $R$ and $M$
  be a free $R$ module of rank $n \ge 2$, and $N=(R \oplus M)$. Then
$\ETrans(N,IN) = \Trans(N,IN) = \E_{n+1}(R,I)$.
\end{lem}

Proof: Let $I_{n+1}$ denote the identity matrix of size $n+1$. Note 
that when $M$ is a free $R$ module, an element of
$\Trans(N, IN)$ looks like $I_{n+1} + v^t w$, for some $v,w \in
R^{n+1}$. Among $v$ and $w$, one of them is unimodular and 
the other belongs to $I^{n+1} ~
(\subseteq R^{n+1})$. Also, note that $\langle v, w \rangle = v w^t =
0$. Therefore $\Trans(N,IN) \subseteq \E_{n+1}(R,I)$ (follows from
\cite{18}, Corollary 1.2 and Lemma 1.3).

Using definition of $\ETrans(N,IN)$, $\E_{n+1}(R,I)$ and the fact that
in the free case elementary transvections $E_x$ and $E^*_\tau$ of
$N$ are of the form $\big( \begin{smallmatrix} 1 & x \\ 0
  & I_n \end{smallmatrix} \big)$, and $\big( \begin{smallmatrix} 1 & 0
  \\ y^t & I_n \end{smallmatrix} \big)$, respectively, we get
$\ETrans(N,IN) \subseteq \E_{n+1}(R,I)$. By (\cite{vdK1}, Lemma 2.2)
and (\cite{SV}, Lemma 2.7(a)) we have $\E_{n+1}(R,I) \subseteq
\ETrans(N,IN)$, hence $\ETrans(N,IN) =\E_{n+1}(R,I)$. We have
$\E_{n+1}(R,I) = \ETrans(N,IN) \subseteq \Trans(N,IN) \subseteq
\E_{n+1}(R,I)$, and hence the result follows.  \hfill{$\square$}

We now establish dilation principle for the elementary transvection
group in the relative case w.r.t. an ideal of the ring. Dilation
principle in the absolute case is proved in (\cite{bbr}, Proposition
3.1).

\begin{lem} \label{linear,dil,rel} Let $I$ be an ideal of $R$ and $M$ 
be a finitely generated module of $R$. Suppose that $a$ is a
non-nilpotent element of $R$ such that $M_a$ is a free $R_a$-module of
rank $n \ge 2$. Let $N=(R \oplus M)$. Let $\alpha(X) \in \Aut(N[X])$
with $\alpha(0)=Id.$, and $\alpha(X)_a \in
\ETrans(N_a[X],IN_a[X])$. Then there exists $\alpha^*(X) \in
\ETrans(N[X],IN[X])$, with $\alpha^*(0) = Id.$, such that
$\alpha^*(X)$ localises to $\alpha(bX)$, for $b \in (a^d)$, $d \gg 0$.
\end{lem}

Proof: Using previous lemma we get
$\ETrans(N_a[X],IN_a[X])=\E_{n+1}(R_a[X],I_a[X])$ and using Lemma
\ref{vanderk1} we get
$\E_{n+1}(R_a[X],I_a[X])=\E_{n+1}^1(R_a[X],I_a[X]) \cap
GL_{n+1}(R_a,I_a)$. Therefore, one can write $\alpha(X)_a = \prod_t
\gamma_t E_{i_t j_t}(X f_t(X)) \gamma_t^{-1}$, where either $i_t=1$ or
$j_t=1$, and $\gamma_t \in E_{n+1}^1(R_a, I_a)$. Note that $f_t(X) \in
R_a[X]$ if $i_t=1$, and $f_t(X) \in I_a[X]$ if $j_t=1$. Using Lemma
\ref{ness4-E1} we write
\begin{eqnarray*}
\alpha(Y^{4^r}X)_a &=& \prod_k E_{i_k j_k}(Y h_k(X,Y)/ a^m),
\end{eqnarray*}
where either $i_k=1$ or $j_k=1$. Note that $h_k(X,Y) \in
R[X,Y]$ if $i_k=1$, and $h_k(X,Y) \in I[X,Y]$ if $j_k=1$.

Since $N_a$ is a free $R_a$ module we have $N_a \cong R_a^{n+1} \cong
N_a^*$. Let $p_1, \ldots, p_{n+1}$ be the standard
basis of $N_a$, $p_1^*, \ldots, p_{n+1}^*$ be the standard basis of
$N_a^*$, and $e_1, \ldots, e_{n+1}$ be the standard basis of
$R_a^{n+1}$. Let $s \ge 0$ be an integer such that $\widetilde{p_i} =
a^s p_i \in N$, and $\widetilde{p_i^*} = a^s p_i^* \in N^*$, for all
$i$.  Note that
\begin{eqnarray*}
\alpha(Y^{4^r} X)_a &=& \prod_k \big( Id. + (Y h_k(X,Y)/ a^m) ~ e_{i_k}^t . ~ e_{j_k} \big) \\
&=& \prod_k \big( Id. + (Y h_k(X,Y)/ a^m) ~ p_{i_k}^* . ~ p_{j_k} \big)
\end{eqnarray*}

Let $m \ge 0$ be the maximum power of $a$ appearing in the
denominators of each $h_k$. Let $d'=2s+m$. Now,
\begin{eqnarray*}
\alpha((a^{d'}Y)^{4^r} X)_a &=& \prod_k \big( Id. + a^{2s}Y
h_k(X,a^{d'}Y)) ~ p_{i_k}^* . ~ p_{j_k} \big) \\ &=& \prod_k \big( Id. +
Y h_k(X,a^{d'}Y)) ~ a^s p_{i_k}^* . ~ a^s p_{j_k} \big).
\end{eqnarray*}

Substituting $Y=1$ we get $\alpha(a^d X)_a = \prod_k \big( Id. + 
h'_k(X) a^s p_{i_k}^* . a^s p_{j_k} \big)$. Let us set $\alpha^*(X) =
\prod_k \big( Id. +  h'_k(X) \widetilde{p_{i_k}^*}
. \widetilde{p_{j_k}} \big)$, where either $i_k=1$ or $j_k=1$. 
Note that $\alpha^*(X)$ belongs to
$\ETrans(N[X]) \cap \Aut(N[X], IN[X]$, where $\Aut(N[X], IN[X] =
ker(\Aut(N[X]) \longrightarrow \Aut(N[X]/IN[X]))$. Generators of
$\ETrans(N[X])$ have the splitting property as $E_{x+y}=E_x E_y$ 
and $E^*_{\tau+\rho}=E^*_{\tau} E^*_{\rho}$. Therefore, using 
argument similar to Lemma \ref{vanderk1} we can show that 
$\alpha^*(X) \in \ETrans(N[X], IN[X])$. It is clear from the construction 
that $\alpha^*(0)=Id.$ and $\alpha^*(X)$ localises to $\alpha(bX)$, 
for $b \in (a^d)$, $d \gg 0$.
\hfill{$\square$}

\begin{lem} \label{linear,LG,rel} Let $I$ be an ideal of $R$ and let
  $M$ be a finitely generated projective $R$-module of rank $n \ge 2$.
  Let $N=(R \oplus M)$. Let $\alpha(X) \in \Aut(N[X])$, with $\alpha(0)
  =Id$. If for each maximal ideal $\gm$ of $R$, $\alpha(X)_\gm \in
  \ETrans(N_\gm[X],IN_\gm[X])$, then $\alpha(X) \in
  \ETrans(N[X],IN[X])$.
\end{lem}

Proof: One can suitably choose an element $a_\gm$ from $R \setminus
\gm$ such that $\alpha(X)_{a_\gm}$ belongs to
$\ETrans(N_{a_\gm}[X],IN_{a_\gm}[X])$. Let us set $\gamma(X,Y) =
\alpha(X+Y) \alpha(Y)^{-1}$. Note that $\gamma(X,Y)_{a_\gm}$ belongs
to $\ETrans(N_{a_\gm}[X,Y],IN_{a_\gm}[X,Y])$, and $\gamma(0,Y) =
Id$. From Lemma \ref{linear,dil,rel} it follows that $\gamma(b_\gm
X,Y) \in \ETrans(N[X,Y],IN[X,Y])$, for $b_\gm \in (a_\gm^d)$, where $d
\gg 0$. Note that the ideal generated by $a_\gm^d$'s is the whole ring
$R$. Therefore, $c_1 a_{\gm_1}^d+ \ldots + c_k a_{\gm_k}^d = 1 $, for
some $c_i \in R$. Let $b_{m_i}= c_i a_{m_i}^d \in (a_{m_i}^d)$. It is
easy to see that $\alpha(X)=\prod_{i=1}^{k-1}\gamma(b_{m_i}X,T)
\gamma(b_{m_k}X,0)$, where $T_i = b_{m_{i+1}}X+ \cdots +
b_{m_k}X$. Each term in the right hand side of this expression belongs
to $\ETrans(N[X], IN[X])$ and hence $\alpha(X) \in \ETrans(N[X],
IN[X])$.  \hfill{$\square$}

\begin{thm} \label{linear,LGaction-rel} Let $I$ be an ideal of $R$ and
  let $M$ be a finitely generated projective $R$-module of rank $n \ge
  2$. Let $N=(R \oplus M)$. Let $q(X) \in \Um(N[X],IN[X])$,
  where $q(X)$ is of the form $(a(X),p(X))$. If for each maximal ideal
  $\gm$ of $R$, $q(X) \in q(0) \ETrans(N_\gm[X],IN_\gm[X])$, then
  $q(X) \in q(0) \ETrans(N[X],IN[X])$.  
\end{thm}

Proof: For each maximal ideal $\gm$ of $R$, we get $\beta_{(\gm)}(X)$
in $\ETrans(N_\gm[X],IN_\gm[X])$ such that $q(X) \beta_{(\gm)}(X) =
q(0)$.  Let us set $\gamma(X,T) = \beta_{(\gm)}(X+T)
\beta_{(\gm)}(X)^{-1}$.  Clearly $\gamma(X,T) \in
\ETrans(N_\gm[X,T],IN_\gm[X,T])$. Since there are only finitely many
denominators involved in the expression of $\gamma(X,T)$, there exists
$a_\gm \in R \setminus \gm$ such that $\gamma(X,T)$ belongs to
$\ETrans(N_{a_\gm}[X,T],IN_{a_\gm}[X,T])$ and $\gamma(X,0)=Id$. Using
Lemma \ref{linear,dil,rel} it follows that $\gamma(X,b_{\gm}T) \in$
$\ETrans(N[X,T],IN[X,T])$ for $b_\gm \in (a_\gm^d)$, $d \gg 0$. We
have $q(X+b_{\gm}T) \gamma(X, b_{\gm}T) = q(X+b_{\gm}T)
\beta_{(\gm)}(X+b_{\gm}T) \beta_{(\gm)}(X)^{-1} = q(0)
\beta_{(\gm)}(X)^{-1} = q(X)$.

Note that the ideal generated by ${a_\gm^d}$'s is the whole ring $R$.
Therefore $c_1 a_{\gm_1}^d + \cdots + c_k a_{\gm_k}^d = 1$, for some
$c_i \in R$. Let $b_{m_i} = c_i a_{m_i}^d \in (a_{m_i}^d)$. In the
above equation replacing $X$ by $b_{\gm_2}X + \cdots + b_{\gm_k}X$ and
replacing $b_m T$ by $b_{m_1} X$ we get, $$q(X) = q (b_{\gm_1}X + b_{\gm_2}X +
\cdots + b_{\gm_k}X) \in q (b_{\gm_2}X + \cdots + b_{\gm_k}X) ~ {\rm
  ETrans}(N[X],IN[X]).$$ Again, in the above equation replacing $X$ by
$b_{\gm_3}X + \cdots + b_{\gm_k}X$ and replacing $b_{\gm}T$ by $b_{\gm_2}X$ we
get, $q( b_{\gm_2}X + \cdots + b_{\gm_k}X) \in q(b_{\gm_3}X + \cdots +
b_{\gm_k}X) {\rm ETrans}(N[X],IN[X])$. Continuing in this way we get
$q(b_{m_k}X + 0) \in q(0) {\rm ETrans}(N[X],IN[X])$. Combining all
these we get $q(X) \prod_{i=1}^{k-1} \gamma(b_{\gm_{i+1}}X + \cdots +
b_{\gm_k}X, b_{\gm_i}X) \gamma(0, b_{\gm_k}X) = q(0)$, where
  $\prod_{i=1}^{k-1} \gamma(b_{\gm_{i+1}}X + \cdots + b_{\gm_k}X,
  b_{\gm_i}X) \gamma(0, b_{\gm_k}X) \in {\rm ETrans}(N[X],IN[X])$.
    \hfill{$\square$}

\medskip
Before we establish equality of the transvection group and elementary
transvection group, we prove a lemma to show that every transvection is
homotopic to identity. Note that for any abelian group $G$, we
understand by $G[X]$ the abelian group of all polynomials in $X$ with
coefficients in $G$.

\begin{lem} \label{lin-trans-hom-to-identity}
Let $M$ be a finitely generated $R$-module and $\alpha \in
\Trans(M)$. Then there exists $\beta(X) \in \Trans(M[X])$ such that
$\beta(1)=\alpha$ and $\beta(0)=Id.$
\end{lem}

Proof: As $\alpha \in \Trans(M)$, it is product of elements of the
form $Id+\pi_q$, where $\pi \in M^*$, $q \in M$ with $\pi(q)=0$. Here
either $\pi$ or $q$ is unimodular. Let $\pi X$ represent $\pi$ {\it
  times} $X$ which belongs to $M^*[X]$ and $q X$ represent $q$ {\it
  times} $X$ which belongs to $M[X]$. Define $\pi X (q) = \pi(q) X$ or
$\pi(qX)=\pi(q)X$. We set $\beta(X)$ to be the product of elements of
the form $Id+\pi_{qX}$ or $Id+\pi X_{q}$, whenever $Id+\pi_q$ appears
in the expression of $\alpha$. This choice depends on whether $\pi$ is
unimodular or $q$ is unimodular. Then $\beta(1)=\alpha$ and
$\beta(0)=Id$.  \hfill{$\square$}


\begin{pr} \label{Trans,equal,rel} Let $I$ be an ideal of $R$. Let $M$ be
  a finitely generated projective $R$-module of rank at least 2, and
  $N=(R \oplus M)$. Then $\Trans(N,IN) = \ETrans(N,IN)$.
\end{pr}

Proof: Note that $\ETrans(N,IN) \subseteq \Trans(N,IN)$. Let us
consider an element $\alpha \in \Trans(N,IN)$. By Lemma
\ref{lin-trans-hom-to-identity} there exists $\alpha(X) \in
\Trans(N[X],IN[X])$ such that $\alpha(1) = \alpha$ and $\alpha(0)=Id$.
Let $\gm$ be a maximal ideal of $R$. We have $\alpha(X)_\gm \in
\Trans(N_\gm[X],IN_\gm[X]) = \ETrans(N_\gm[X],IN_\gm[X])$ (see Lemma
\ref{linear,free,rel}). This is true for all maximal ideal $\gm$ of
$R$ and hence $\alpha(X) \in \ETrans(N[X],IN[X])$ by Lemma
\ref{linear,LG,rel}. Substituting $X=1$ we get $\alpha \in
\ETrans(N,IN)$. Therefore $\Trans(N,IN) \subseteq \ETrans(N,IN)$.
\hfill{$\square$}

\section{\large{Symplectic Modules and Symplectic Transvections}}

\begin{de} {\bf Alternating Matrix:} A matrix in ${\rm M}_n (R)$ is said
  to be {\it alternating} if it has the form $\nu - \nu^t$, where $\nu
  \in {\rm M}_n(R)$. \tn{(It follows that its diagonal elements are
  zeros.)}
\end{de}

\begin{lem} \label{local,rel} Let $(R,\gm)$ be a local ring and $I$ be an
  ideal of $R$. Let $\varphi$ be an alternating matrix of Pfaffian 1
  over $R$, and $\varphi \equiv \psi_n ~({\rm mod}~I)$. Then $\varphi$
  is of the form $$(1 \perp \varepsilon)^t \psi_n (1 \perp
  \varepsilon),$$ for some $\varepsilon \in {\rm E}_{2n-1}(R,I)$.
\end{lem}

Proof: We will prove the result using induction on $n$. When $\varphi$
is of size $2 \times 2$, the result is true. Let us assume that the
result is true for alternating matrix of size $2(n-1) \times 2(n-1)$,
i.e, for any $2(n-1) \times 2(n-1)$ alternating matrix $\varphi^*$ of
Pfaffian 1 with $\varphi^* \equiv \psi_{n-1} ~({\rm mod}~I)$, we have
$\eta \in \E_{2n-3}(R,I)$ such that $\varphi^* = (1 \perp \eta)^t ~
\psi_{n-1} ~ (1 \perp \eta)$.

We will prove the result for alternating matrix
$\varphi$ of size $2n \times 2n$.  Let
$\varphi = \big( \begin{smallmatrix} 0 & a \\ -a^t & \alpha
\end{smallmatrix} \big) \equiv  \psi_n ({\rm mod}~I)$,
where $a \in$ Um$_{2n-1}(R,I)$ and $\alpha$ is alternating matrix of
size $(2n-1) \times (2n-1)$. Note that
$\alpha \equiv \big( \begin{smallmatrix} 0 & 0 \\ 0 & \psi_{n-1} 
\end{smallmatrix} \big) ({\rm mod}~I)$.

As $R$ is local ring we have $a = e_1 \beta$, where $\beta \in$
E$_{2n-1}(R,I)$ (see Remark \ref{unimodularOverLocalRing}). Hence $(1
\perp \beta^t)^{-1} \varphi (1 \perp \beta)^{-1} =
\big( \begin{smallmatrix} 0 & e_1 \\ -e_1^{t} &
  \gamma \end{smallmatrix} \big)$, where $\gamma = (\beta^t)^{-1} ~
\alpha ~ \beta^{-1}$. Note that $\gamma$ is an alternating
matrix. Therefore $\gamma$ can be written as
$\big( \begin{smallmatrix} 0 & b \\ -b^{t} & \varphi^*
\end{smallmatrix} \big)$. Note that
$\overline{\gamma} = (\overline{\beta}^t)^{-1} \overline{\alpha}
\overline{\beta}^{-1} \equiv \big( \begin{smallmatrix} 0 & 0 \\ 0 &
  \psi_{n-1} \end{smallmatrix} \big) ({\rm mod}~I)$, and hence $b \in
I^{2n-2}$ and $\varphi^* \equiv \psi_{n-1} ~ ({\rm mod}~I)$. Now
$\Big( \begin{smallmatrix} 1 & 0 & 0 \\ 0 & 1 & -b {\varphi^*}^{-1}
  \\ 0 & 0 & I_{2n-2} \end{smallmatrix} \Big) (1 \perp \beta^t)^{-1}
\varphi (1 \perp \beta)^{-1} \Big( \begin{smallmatrix} 1 & 0 & 0 \\ 0
  & 1 & -b {\varphi^*}^{-1} \\ 0 & 0 & I_{2n-2} \end{smallmatrix}
\Big)^t = \Big( \begin{smallmatrix} 0 & 1 & 0 \\ -1 & 0 & 0 \\ 0 & 0 &
  \varphi^*
\end{smallmatrix} \Big)$.

Let us call the matrix $((I_3 \perp \eta)^{-1})^t
\Big( \begin{smallmatrix} 1 & 0 & 0 \\ 0 & 1 & -b {\varphi^*}^{-1}
  \\ 0 & 0 & I_{2n-2} \end{smallmatrix} \Big) (1 \perp \beta^t)^{-1} =
((1 \perp \varepsilon)^{-1})^t$. Note that $\varepsilon \in
\E_{2n-1}(R,I)$. Using induction hypothesis we get $((1 \perp
\varepsilon)^{-1})^t \varphi (1 \perp \varepsilon)^{-1} = ((I_3 \perp
\eta)^{-1})^t \big( \begin{smallmatrix} 0 & 1 & 0 \\ -1 & 0 & 0 \\ 0 &
  0 & \varphi^*
\end{smallmatrix} \big) (I_3 \perp \eta)^{-1} = \psi_n$,
and hence $\varphi = (1 \perp \varepsilon)^t \psi_n (1 \perp
\varepsilon)$.  Therefore the result is established.
\hfill{$\square$}

\begin{rmk} \label{alternating} The condition that the alternating
  matrices are of Pfaffian one can be replaced by the weaker condition
  that the alternating matrix be invertible, and congruent to
  $(u \psi_1 \perp \psi_{n-1}) ~({\rm mod}~I)$, where $u$ = Pfaffian
  $\varphi$, is of the form $(1 \perp E)^t (u \psi_1 \perp \psi_{n-1})
  (1 \perp E)$, for some relative elementary matrix $E$. 
\end{rmk}

\begin{rmk} \label{local-rel} Let $\varphi$ be an alternating matrix
  of Pfaffian $1$, over $R$. Let $\gm$ be a maximal ideal of $R$ and
  $R_\gm$ be the local ring at $\gm$. We will get $\varepsilon(\gm)
  \in {\rm E}_{2n-1}(R_\gm)$ such that over $R_\gm$ we have $\varphi =
  (1 \perp \varepsilon(\gm))^t \psi_n (1 \perp \varepsilon(\gm))$. Let
  $a$ be the product of denominators of all the entries of
  $\varepsilon(\gm)$. Clearly $a$ is not in $\gm$.  Hence we get
  $\varepsilon$ from ${\rm E}_{2n-1}(R_a)$ such that $\varphi = (1
  \perp \varepsilon)^t \psi_n (1 \perp \varepsilon)$.

{\it Note:} {\tn When dealing with relative case w.r.t. an ideal $I$ of $R$,
  we will always assume that the alternating matrix $\varphi$ of
  Pfaffian 1 is congruent to $\psi_n ~({\rm mod} ~ I)$. Using Lemma
  \ref{local,rel} and arguing as above we get that over ring $R_\gm$
  we have $\varphi = (1 \perp \varepsilon)^t \psi_n (1 \perp
  \varepsilon)$, where $\varepsilon \in {\rm E}_{2n-1}(R_a,I_a)$, for
  some $a \notin \gm$. We will constantly use this fact without even
  referring to it!}
\end{rmk}

\begin{de} 
{\rm 
A {\bf symplectic $R$-module} is a pair $(P,\langle , \rangle)$,
where $P$ is a finitely generated projective $R$-module of even rank
and $\langle , \rangle: P \times P \lra R$ is a non-degenerate (i.e, $P
\cong P^*$ by $x \lra \langle x , -\rangle$) alternating bilinear form.
}
\end{de}

\begin{de}
{\rm
Let $(P_1,\langle , \rangle_1)$ and $(P_2,\langle , \rangle_2)$ be
two  symplectic $R$-modules. Their {\bf orthogonal sum} is the pair
$(P,\langle , \rangle)$, where $P=P_1 \oplus P_2$ and the inner product
is defined by $\langle (p_1,p_2),(q_1,q_2)\rangle = \langle
p_1,q_1\rangle_1 + \langle p_2,q_2\rangle_2$.

There is a non-degenerate bilinear form $\langle , \rangle$ on the
$R$-module $\mathbb{H}(R)= R \oplus R^*$, namely $\langle (a_1,f_1),
(a_2,f_2) \rangle = f_2(a_1) - f_1(a_2)$.  }
\end{de}

\begin{de}
{\rm An {\bf isometry} of a symplectic module $(P,\langle , \rangle)$ is
  an automorphism of $P$ which fixes the bilinear form. The group of
  isometries of $(P, \langle , \rangle)$ is denoted by
  $\Sp(P, \langle , \rangle)$.  }
\end{de}

\begin{de} 
{\rm 
In \cite{bass2} Bass has defined a {\bf symplectic transvection}
of a symplectic module $P$ to be an automorphism of the form
\begin{eqnarray*}
\sigma(p) &=& p + \langle u , p \rangle v + \langle v , p \rangle u + \alpha
\langle u , p \rangle u,
\end{eqnarray*}
where $\alpha \in R$, $u,v \in P$ are fixed elements with $\langle
u,v \rangle=0$, and either $u$ or $v$ is unimodular. It is easy to
check that $\langle \sigma(p), \sigma(q) \rangle = \langle p, q
\rangle$ and $\sigma$ has an inverse $\tau(p) = p - \langle u, p
\rangle v - \langle v, p \rangle u - \alpha \langle u, p \rangle u$.

The subgroup of $\Sp(P,\langle , \rangle)$ generated by the symplectic
transvections is denoted by $\Trans_{\Sp}(P, \langle , \rangle)$ (see
\cite{swan}, Page 35).  }
\end{de}

{\bf Now onwards $Q$ will denote $(R^2 \oplus P)$ with induced form
  on $(\mathbb{H}(R)~\oplus~P)$, and $Q[X]$ will denote $(R[X]^2 \oplus
  P[X])$ with induced form on $(\mathbb{H}(R[X])~\oplus~P[X])$.}

\begin{de} {\rm The symplectic transvections of $Q=(R^2 \oplus P)$ of 
the form
    \begin{eqnarray*}
      (a, b, p) & \mapsto & (a, b - \langle p, q \rangle - \alpha a, p-aq),
     \end{eqnarray*} 
or of the form 
    \begin{eqnarray*}
      (a, b, p) & \mapsto & (a + \langle p, q \rangle + \beta b, b,  p-bq), 
    \end{eqnarray*} 
where $\alpha, \beta \in R$ and $q \in P$, 
are called {\bf elementary symplectic transvections}. Let us denote the first 
isometry by $\rho(q, \alpha)$ and the second one by $\mu(q, \beta)$. It can 
be verified that the elementary symplectic transvections are symplectic 
transvections on $Q$. Indeed, consider $(u, v) =((0,1,0),(0,0,q))$  to get 
$\rho(q, \alpha)$ and consider $(u, v)= ((-1,0,0),(0,0,q))$ to get $\mu(q, \beta)$.

The subgroup of $\Trans_{\Sp}(Q, \langle , \rangle)$ generated by
elementary symplectic transvections is denoted by $\ETrans_{\Sp}(Q,
\langle , \rangle)$.}
\end{de}

\begin{de} 
{\rm Let $I$ be an ideal of $R$. The group of {\bf relative symplectic
    transvections} w.r.t. an ideal $I$ is generated by the symplectic
  transvecions of the form $\sigma(p) = p + \langle u, p \rangle v +
  \langle v, p \rangle u + \alpha \langle u, p \rangle u$, where $\alpha
  \in I$ and $u \in P$, $v \in IP$ are fixed elements with $\langle
  u, v \rangle=0$. The group generated by relative symplectic
  transvections is denoted by $\Trans_{\Sp}(P,IP, \langle , \rangle)$.}
\end{de}

\begin{de} {\rm Let $I$ be an ideal of $R$. The elementary symplectic
transactions of $Q$ of the form $\rho(q, \alpha), \mu(q, \beta)$, where
$q \in IP$ and $\alpha, \beta \in I$ are called {\bf relative elementary 
symplectic transvections} w.r.t. an ideal $I$.

The subgroup of $\ETrans_{\Sp}(Q, \langle , \rangle)$ generated by
relative elementary symplectic transvections is denoted by
$\ETrans_{\Sp}(IQ,\langle , \rangle )$. The normal closure of
$\ETrans_{\Sp}(IQ,\langle , \rangle )$ in
$\ETrans_{\Sp}(Q, \langle , \rangle)$ is denoted by
$\ETrans_{\Sp}(Q,IQ, \langle , \rangle)$ }.
\end{de}

\begin{rmk} \label{free} Let $P=\oplus_{i=1}^{2n} Re_i$ be a free $R$-module. 
The non-degenerate alternating bilinear form $\langle,\rangle$ on $P$
corresponds to an alternating matrix $\varphi$ with Praffian 1 with
respect to the basis $\{e_1, e_2, \ldots, e_{2n} \}$ of $P$ and we write
$\langle p, q \rangle = p \varphi q^t$.

  In this case the symplectic transvection $\sigma(p) = p + \langle u,
  p \rangle v + \langle v, p \rangle u + \alpha \langle u, p \rangle
  u$ corresponds to the matrix $(I_{2n} - \varphi u^t v - \varphi v^t
  u)(I_{2n} - \alpha \varphi u^t u)$ and the group generated by them
  is denoted by $\Trans_{\Sp}(P, \langle , \rangle_{\varphi})$.

Also in this case $\ETrans_{\Sp}(Q, \langle , \rangle_{\psi_1 \perp \varphi})$ will
be generated by the matrices of the form $\rho_{\varphi}(q, \alpha) =
\Big( \begin{smallmatrix} 1 & -\alpha & -q \\ 0 & 1 & 0 \\ 0 &
  - \varphi q^t & I_{2n} \end{smallmatrix} \Big)$, and $\mu_{\varphi}(q, \beta) =
\Big( \begin{smallmatrix} 1 & 0 & 0 \\ \beta & 1 & -q \\ \varphi q^t&
0 & I_{2n} \end{smallmatrix} \Big)$.

Note that for $q=(q_1, \ldots, q_{2n}) \in R^{2n}$, and for the
standard alternating matrix $\psi_n$, we have
\begin{eqnarray} \label{relatn5}
\rho_{\psi_n}(q, \alpha) &=& se_{12}(-\alpha) \prod_{i=3}^{2n+2}
se_{1i}(-q_{i-2}), \\ 
\label{relatn6}
\mu_{\psi_n}(q, \beta) &=& se_{21}(\beta)
\prod_{i=3}^{2n+2} se_{i1}((-1)^{i+1} q_{\sigma(i-2)}).  
\end{eqnarray}
\end{rmk}

{\it  We shall implicitly use the assumptions and notations in the
  statement of Remark \ref{free} in the sequel below.}

\begin{lem} \label{kopeiko} 
Let $R$ be a commutative ring and $I$ be an ideal of $R$. Let $P$
be a free $R$-module of rank $2n$, $n \ge 1$. If $\varphi = \psi_n$,
the standard alternating matrix, then $\Trans_{\Sp}(Q,IQ, \langle ,
\rangle_{\psi_{n+1}}) = \ESp_{2n+2}(R,I)$.
\end{lem}

Proof: For proof see (\cite{29}, Lemma 1.10).

\begin{lem} \label{free,psi} Let $R$ be a commutative ring with
  $R=2R$, and let $I$ be an ideal of $R$. Let $P$ be a free $R$-module
  of rank $2n$, $n \ge 1$. If $\varphi = \psi_n$, the standard
  alternating matrix, then $\ETrans_{\Sp}(Q,IQ,
  \langle , \rangle_{\psi_{n+1}}) = \ESp_{2n+2}(R,I)$.
\end{lem}

Proof: We first show $\ETrans_{\Sp}(Q,IQ, \langle ,
\rangle_{\psi_{n+1}})$ is a subset of $\ESp_{2n+2}(R,I)$. It is enough
to show that an element of the form $s_{\psi_n}(q_1, \alpha_1)
t_{\psi_n}(q_2, \beta) s_{\psi_n}(q_1, \alpha)^{-1}$ is in
$\ESp_{2n}(R,I)$, for $q_1 \in R^{2n}$, $q_2 \in I^{2n} (\subseteq
R^{2n})$, $\alpha \in R$, and $\beta \in I$. Here $t_{\psi_n}$ and
$s_{\psi_n}$ represents either $\rho_{\psi_n}$ or
$\mu_{\psi_n}$. Using equations (\ref{relatn5}), (\ref{relatn6}) we
get $\ETrans_{\Sp}(Q,IQ, \langle , \rangle_{\psi_{n+1}}) \subseteq
\ESp_{2n}(R,I)$.

To show the other inclusion we recall the equivalent definition of the
relative group which says that $\ESp_{2n}(R,I)$ is the smallest normal
subgroup of $\ESp_{2n}(R)$ containing $se_{21}(x)$, where $x \in I$.
We have $g se_{21}(x) g^{-1} \in
\ETrans_{\Sp}(Q,IQ,\langle , \rangle_{\psi_{n+1}})$, hence
$\ESp_{2n}(R,I) \subseteq \ETrans_{\Sp}(Q,IQ, \langle ,
\rangle_{\psi_{n+1}})$. Therefore, the equality is established.
\hfill{$\square$}

\begin{lem} \label{phi=phi*,ab}
Let $P$ be a free $R$-module of rank $2n$. Let $(P,\langle ,
\rangle_{\varphi})$ and $(P,\langle , \rangle_{\varphi^*})$ be two
symplectic $R$-modules with $\varphi= (1 \perp \varepsilon)^t ~
\varphi^* ~ (1 \perp \varepsilon)$, for some $\varepsilon \in
\E_{2n-1}(R)$. Then
\begin{eqnarray*} 
\Trans_{\Sp}(P,\langle , \rangle_{\varphi}) &=& (1 \perp
\varepsilon)^t ~ \Trans_{\Sp}(P,\langle , \rangle_{\varphi^*}) ~ ((1
\perp \varepsilon)^t)^{-1},\\ 
\ETrans_{\Sp}(Q, \langle , \rangle_{\psi_1 \perp \varphi}) &=& (I_3 \perp
\varepsilon)^t ~ \ETrans_{\Sp}(Q,\langle , \rangle_{\psi_1 \perp \varphi^*}) ~ ((I_3
\perp \varepsilon)^t)^{-1}.
\end{eqnarray*}
\end{lem}

Proof: In the free case for symplectic transvections we have
\begin{eqnarray*}
 & & (I_{2n} - \varphi u^t v - \varphi v^t u )(I_{2n} - \alpha \varphi
  u^t u)\\ & = & (1 \perp \varepsilon)^t ~ (I_{2n} - \varphi^*
  \tilde{u}^t \tilde{v} - \varphi^* \tilde{v}^t \tilde{u}) ~ (I_{2n} -
  \alpha \varphi^* \tilde{u}^t \tilde{u}) ~ ((1 \perp \varepsilon)^t)^{-1},
\end{eqnarray*}
where $\tilde{u}=u(1 \perp \varepsilon)^t$ and $\tilde{v}=v(1 \perp
\varepsilon)^t$. Hence the first equality follows. 

For elementary symplectic transvections we have 
\begin{eqnarray*}
(I_2 \perp (1 \perp \varepsilon))^t \rho_{\varphi^*}(q, \alpha) ((I_2
  \perp (1 \perp \varepsilon))^t)^{-1} = \rho_{\varphi} (q ((1 \perp
  \varepsilon)^t)^{-1}, \alpha),\\ (I_2 \perp (1 \perp
  \varepsilon))^t \mu_{\varphi^*}(q, \beta) ((I_2 \perp (1 \perp
  \varepsilon))^t)^{-1} = \mu_{\varphi} (q ((1 \perp
  \varepsilon)^t)^{-1}, \beta),
\end{eqnarray*}
hence the second equality follows.
\hfill{$\square$}

\begin{lem} \label{phi=phi*,rel} Let $I$ be an ideal of $R$ and $P$ be
  a free $R$-module of rank $2n$. Let $(P,\langle , \rangle_{\varphi})$
  and $(P,\langle , \rangle_{\varphi^*})$ be two symplectic $R$-modules
  with $\varphi= (1 \perp \varepsilon)^t ~ \varphi^* ~ (1 \perp
  \varepsilon)$, for some $\varepsilon \in \E_{2n-1}(R,I)$. Then
\begin{eqnarray*} 
\Trans_{\Sp}(P,IP,\langle , \rangle_{\varphi}) &=& (1 \perp
\varepsilon)^t ~ \Trans_{\Sp}(P,IP,\langle , \rangle_{\varphi^*}) ~ ((1
\perp \varepsilon)^t)^{-1},\\ 
\ETrans_{\Sp}(Q,IQ,\langle , \rangle_{\psi_1 \perp \varphi}) &=& (I_3 \perp
\varepsilon)^t ~ \ETrans_{\Sp}(Q,IQ,\langle , \rangle_{\psi_1 \perp \varphi^*}) ~ ((I_3
\perp \varepsilon)^t)^{-1}.
\end{eqnarray*}
\end{lem}

Proof: Using the three equations appearing in the proof of Lemma
\ref{phi=phi*,ab}, we get these equalities.
\hfill{$\square$}

\begin{pr} \label{local-case,rel} Let $R$ be a commutative ring with
  $R=2R$, and let $I$ be an ideal of $R$. Let $(P,\langle ,
  \rangle_{\varphi})$ be a symplectic $R$-module with $P$ free of rank
  $2n$, $n \ge 1$. Let $\varphi=(1 \perp \varepsilon)^t ~ \psi_n ~ (1
  \perp \varepsilon)$, for some $\varepsilon \in \E_{2n-1}(R,I)$. Then
  $\Trans_{\Sp}(Q,IQ,\langle , \rangle_{\psi_1 \perp \varphi}) =
  \ETrans_{\Sp}(Q,IQ,\langle , \rangle_{\psi_1 \perp \varphi})$.
\end{pr}

Proof: Using Lemma \ref{kopeiko}, Lemma \ref{free,psi}, and Lemma
\ref{phi=phi*,rel} we get,

\begin{eqnarray*}
\Trans_{\Sp}(Q,IQ,\langle , \rangle_{\psi_1 \perp \varphi}) &=& (I_3
\perp \varepsilon)^t~ \Trans_{\Sp}(Q,IQ,\langle , \rangle_{\psi_{n+1}})
  ~ ((I_3 \perp \varepsilon)^t)^{-1}\\ &=& (I_3 \perp \varepsilon)^t ~
  \ESp_{2+2n}(R,I) ~ ((I_3 \perp \varepsilon)^t)^{-1},
\end{eqnarray*}
and
\begin{eqnarray*}
\ETrans_{\Sp}(Q,IQ,\langle , \rangle_{\psi_1 \perp \varphi}) &=& (I_3 \perp
\varepsilon)^t ~ \ETrans_{\Sp}(Q,IQ,\langle , \rangle_{\psi_{n+1}}) ~ ((I_3
\perp \varepsilon)^t)^{-1}\\ &=& (I_3 \perp \varepsilon)^t ~ \ESp_{2+2n}(R,I)
~ ((I_3 \perp \varepsilon)^t)^{-1},
\end{eqnarray*}
and hence the equality is established. 
\hfill{$\square$}

\begin{rmk} \label{localcase} In view of above two lemmas, for any
    symplectic module $(P,\langle , \rangle_{\varphi})$ over a local
    ring $(R,\gm)$, we have $\Trans_{\Sp}(Q,IQ,\langle,\rangle_{\psi_1
      \perp \varphi} ) = \ETrans_{\Sp}(Q,IQ, \langle , \rangle_{\psi_1
      \perp \varphi} )$. Here $I$ is an ideal of the ring $R$.
\end{rmk}

Next we establish dilation principle for elementary symplectic
transvection group.

\begin{lem} \label{rel-dilatn-ETrans} Let $R$ be a commutative ring
  with $R=2R$, and let $I$ be an ideal of $R$. Let $(P,\langle ,
  \rangle)$ be a symplectic $R$-module with $P$ finitely generated
  projective $R$-module of rank $2n$, $n \ge 1$. Suppose that $a$ is a
  non-nilpotent element of $R$ such that $P_a$ is a free $R_a$ module
  and the bilinear form $\langle, \rangle$ corresponds to the
  alternating matrix $\varphi$ (w.r.t. some basis). Assume that
  $\varphi = (1 \perp \varepsilon)^t ~ \psi_n ~ (1 \perp
  \varepsilon)$, for some $\varepsilon \in \E_{2n-1}(R_a,I_a)$. Let
  $\alpha(X) \in \Aut(Q[X])$, with $\alpha(0) = Id$, and $\alpha(X)_a 
  \in \ETrans_{\Sp}(Q_a[X],IQ_a[X],\langle , \rangle_{\psi_1 \perp
    \varphi})$. Then, there exists $\alpha^*(X) \in 
    \ETrans_{\Sp}(Q[X],IQ[X],\langle , \rangle)$, with $\alpha^*(0) =
  Id.$, such that $\alpha^*(X)$ localises to $\alpha(bX)$, for $b \in
  (a^d)$, $d \gg 0$.
\end{lem}

Proof: We have $P_a \cong R_a^{2n}$. Let $e_1, \ldots, e_{2n+2}$ be the standard basis of $Q_a$ with respect to which the bilinear form on $Q_a$ will correspond to $\psi_1 \perp \varphi$. Given that $\alpha(X)_a \in \ETrans_{\Sp}(Q_a[X], IQ_a[X], \langle, \rangle_{\psi_1 \perp \varphi})$ and $\ETrans_{\Sp}(Q_a[X], IQ_a[X], \langle, \rangle_{\psi_1 \perp \varphi}) = (I_3 \perp \varepsilon)^t \ESp_{2n+2}(R_a[X], I_a[X]) ((I_3 \perp \varepsilon)^t)^{-1}$ by Lemma \ref{free,psi}, and Lemma
\ref{phi=phi*,rel}. Therefore, $\alpha(X)_a = (I_3 \perp \varepsilon)^t ~\beta(X)~ ((I_3 \perp \varepsilon)^t)^{-1}$, for some $\beta(X) \in \ESp_{2n+2}(R_a[X], I_a[X])$, with $\beta(0)=Id.$ Note that $\ESp_{2n+2}(R_a[X], I_a[X]) = \ESp_{2n+2}^1(R_a[X], I_a[X]) ~\cap \\ \Sp_{2n+2}(R_a[X], I_a[X])$. Hence we can write $\beta(X) = \prod_t \gamma_t se_{i_t j_t}(X f_t(X)) \gamma_t^{-1}$, where either $i_t =1$, or $j_t=1$, and $\gamma_t \in \ESp_{2n+2}^1(R_a, I_a)$. Note that $f_t(X) \in R_a[X]$, when $i_t=1$ and $f_t(X) \in I_a[X]$, when $j_t=1$. Using Lemma \ref{ness4-E1} we get $\beta(Y^{4^r}X) = \prod_k se_{i_k j_k}(Y h_k(X,Y)/ a^m)$, with either $i_k =1$ or $j_k=1$. Note that $h_k(X,Y) \in R[X,Y]$, when $i_k=1$ and $h_k(X,Y) \in I[X,Y]$, when $j_k=1$. Note that  
\begin{eqnarray*}
se_{12}(Y h_k(X,Y)/ a^m) &=& I_{2n+2} + (Y h_k(X,Y)/ a^m) \psi_{n+1}~ e_2^t ~e_2, \\
se_{1 j_k}(Y h_k(X,Y)/ a^m) &=& I_{2n+2} + (Y h_k(X,Y)/ a^m) \psi_{n+1}~ e_2^t ~e_{j_k} \\
&& + (Y h_k(X,Y)/ a^m) \psi_{n+1}~ {e_{j_k}}^t ~e_2, ~ {\rm for} ~ j_k \ge 3, \\
se_{21}(Y h_k(X,Y)/ a^m) &=& I_{2n+2} - (Y h_k(X,Y)/ a^m) \psi_{n+1}~ e_1^t ~e_1, \\
se_{i_k 1}(Y h_k(X,Y)/ a^m) &=& I_{2n+2} + ((-1)^{\sigma(i_k)}Y h_k(X,Y)/ a^m) \psi_{n+1}~ {e_{\sigma(i_k)}}^t ~e_1 \\
&& + ((-1)^{\sigma(i_k)}Y h_k(X,Y)/ a^m) \psi_{n+1}~ e_1^t ~e_ {\sigma(i_k)}, ~ {\rm for} ~ i_k \ge 3. \\
\end{eqnarray*}

Let $\varepsilon_1, \ldots, \varepsilon_{2n}$ be the rows of the matrix $((1 \perp \varepsilon)^t)^{-1} \in \E_{2n}(R_a, I_a)$. Let $\widetilde{e_i}$ denote the row vector $e_i ((I_3 \perp \varepsilon)^t)^{-1}$ of length $2n+2$. Note that $\widetilde{e_1}=e_1, \widetilde{e_2}=e_2$, and $\widetilde{e_i} = (0, 0, \varepsilon_{i-2})$, for $i \ge 3$. Using Lemma \ref{phi=phi*,rel}  we can write $\alpha(Y^{4^r}X)_a$ as product of elements of the form
\begin{eqnarray*}
&I_{2n+2} + (Y h_k(X,Y)/ a^m) \left ( \begin{smallmatrix} \psi_1 & 0 \\ 0 & \varphi \end{smallmatrix} \right ) \widetilde{e_2}^t \widetilde{e_2} = \rho_{\varphi}(0, Y h_k(X,Y)/ a^m), &\\
& I_{2n+2} + (Y h_k(X,Y)/ a^m) \left ( \begin{smallmatrix} \psi_1 & 0 \\ 0 & \varphi \end{smallmatrix} \right ) \widetilde{e_2}^t \widetilde{e_{j_k}} + (Y h_k(X,Y)/ a^m) \left ( \begin{smallmatrix} \psi_1 & 0 \\ 0 & \varphi \end{smallmatrix} \right ) \widetilde{e_{j_k}}^t \widetilde{e_2} & \\
&= \rho_{\varphi}((Y h_k(X,Y)/ a^m) \varepsilon_{j_k -2}, 0), & \\
&I_{2n+2} - (Y h_k(X,Y)/ a^m) \left ( \begin{smallmatrix} \psi_1 & 0 \\ 0 & \varphi \end{smallmatrix} \right ) \widetilde{e_1}^t \widetilde{e_1} = \mu_{\varphi}(0, Y h_k(X,Y)/ a^m), &\\
& I_{2n+2} - (-1)^{\sigma(i_k)} (Y h_k(X,Y)/ a^m) \left ( \begin{smallmatrix} \psi_1 & 0 \\ 0 & \varphi \end{smallmatrix} \right ) \widetilde{e_{\sigma(i_k)}}^t \widetilde{e_1} & \\
&- (-1)^{\sigma(i_k)} (Y h_k(X,Y)/ a^m) \left ( \begin{smallmatrix} \psi_1 & 0 \\ 0 & \varphi \end{smallmatrix} \right ) \widetilde{e_1}^t \widetilde{e_{\sigma(i_k)}} = \mu_{\varphi}((-1)^{\sigma(i_k)}(Y h_k(X,Y)/ a^m) \varepsilon_{\sigma(i_k)-2}, 0), & 
\end{eqnarray*}
 for $i_k, j_k \ge 3$. Let $s \ge 0$ be an integer such that $\widetilde{\varepsilon_i} = a^s \varepsilon_i \in P$ for all $i = 1, \ldots, 2n$. Let $d' = s+m$. Therefore $\alpha((a^{d'}Y)^{4^r} X)_a$ is product of elements of the form
 \begin{eqnarray*}
 & \rho_{\varphi}(0, a^{d'} Y h_k(X, a^{d'} Y)/ a^m), \rho_{\varphi}((a^{d'} Y h_k(X, a^{d'} Y)/ a^m) \varepsilon_{j_k -2}, 0), & \\
 & \mu_{\varphi}(0, a^{d'} Y h_k(X, a^{d'} Y)/ a^m), ~{\rm and}~ \mu_{\varphi}((-1)^{\sigma(i_k)}(a^{d'} Y h_k(X, a^{d'} Y)/ a^m) \varepsilon_{\sigma(i_k)-2}). & 
 \end{eqnarray*}
 
Substituting $Y=1$ we get $\alpha(a^dX)_a$ is product of elements of the forms 
\begin{eqnarray*}
& \rho_{\varphi}(0, a^s h'_k(X)), \rho_{\varphi}(a^s h'_k(X) \varepsilon_{j_k-2}, 0), & \\
& \mu_{\varphi}(0, a^s h'_k(X)), ~{\rm and}~ \mu_{\varphi}((-1)^{\sigma(i_k)} a^s h'_k(X) \varepsilon_{\sigma(i_k)-2}, 0). &
\end{eqnarray*}

Let us set $\alpha^*(X)$ to be the product of elements of the forms 
\begin{eqnarray*}
& \rho (0, a^s h'_k(X)), \rho (h'_k(X) \widetilde{\varepsilon_{j_k-2}}, 0), & \\
& \mu (0, a^s h'_k(X)), ~{\rm and}~  \mu ((-1)^{\sigma(i_k)} h'_k(X) \widetilde{\varepsilon_{\sigma(i_k)-2}}, 0). &
\end{eqnarray*}

Note that $\alpha^*(X)$ belongs to $\ETrans_{\Sp}(Q[X], \langle, \rangle) \cap \Aut(Q[X], IQ[X])$. Generators of $\ETrans_{\Sp}(Q[X])$ have the splitting property as 
\begin{eqnarray*}
\rho(q_1/2, \alpha_1/2) ~ \rho(q_2, \alpha_2) ~ \rho(q_1/2, \alpha_1/2) &=& \rho(q_1 + q_2, \alpha_1 + \alpha_2),\\
\mu(q_1/2, \alpha_1/2) ~ \mu(q_2, \alpha_2) ~ \mu(q_1/2, \alpha_1/2) &=& \mu(q_1 + q_2, \alpha_1 + \alpha_2).
\end{eqnarray*}
Therefore, using argument similar to Lemma \ref{vanderk1} we can show that $\alpha^*(X)$ is in $\ETrans_{\Sp}(Q[X], IQ[X], \langle, \rangle)$. It is clear from construction that $\alpha^*(0)=Id.$ and  $\alpha^*(X)$ localises to $\alpha(bX)$, for some $b \in (a^d), d \gg 0$.
\hfill{$\square$}

\begin{lem} \label{LG-ETrans-rel} Let $R$ be a commutative ring with
  $R=2R$, and let $I$ be an ideal of $R$. Let $(P,\langle , \rangle)$
  be a symplectic $R$-module with $P$ finitely generated projective
  module of rank $2n$, $n \ge 1$. Let $\alpha(X) \in {\Sp}(Q[X],
  \langle , \rangle)$, with $\alpha(0) = Id$.  If for each maximal
  ideal $\gm$ of $R$, $\alpha(X)_\gm \in
  \ETrans_{\Sp}(Q_\gm[X],IQ_\gm[X],\langle , \rangle_{\psi_1 \perp
    \varphi_\gm})$, where $\varphi_\gm \equiv \psi_n ~ ({\rm mod}~I)$,
  then $\alpha(X) \in \ETrans_{\Sp}(Q[X],IQ[X],\langle , \rangle)$.
\end{lem}

Proof: Proof follows arguing similarly as in Lemma \ref{linear,LG,rel}.
\hfill{$\square$}

\begin{thm} \label{LG-ETrans-action-rel} Let $R$ be a commutative ring
  with $R=2R$, and let $I$ be an ideal of $R$. Let
  $(P,\langle , \rangle)$ be a symplectic $R$-module with $P$ a finitely
  generated projective module of rank $2n$, $n \ge 1$. Let
  $q(X)=(a(X),b(X),p(X)) \in \Um(Q[X],IQ[X])$. If for each maximal
  ideal $\gm$ of $R$, we have $q(X) \in q(0)
  ~\ETrans_{\Sp}(Q_\gm[X],IQ_\gm[X],\langle,\rangle_{\psi_1 \perp \varphi_\gm})$,
  where $\varphi_\gm \equiv \psi_n ~ ({\rm mod}~I)$, then $q(X) \in
  q(0) ~\ETrans_{\Sp}(Q[X],IQ[X],\langle , \rangle)$.  
\end{thm}

Proof: Proof follows arguing similarly as in Theorem
\ref{linear,LGaction-rel}.  \hfill{$\square$}

\medskip
We now establish equality of the symplectic transvection group and the
elementary symplectic transvection group. Before that we prove a lemma
to show that symplectic transvections are homotopic to identity.

\begin{lem} \label{symp-trans-hom-to-identity}
Let $(P, \langle, \rangle)$ be a symplectic $R$-module and $\alpha \in
\Trans_{\Sp}(P, \langle, \rangle)$. Then there exists $\beta(X) \in
\Trans_{\Sp}(P[X], \langle, \rangle)$ such that $\beta(1)=\alpha$ and
$\beta(0)=Id.$
\end{lem}

Proof: As $\alpha \in \Trans_{\Sp}(P, \langle, \rangle)$, it is
product of symplectic transvections of the form $\sigma$, where
$\sigma$ takes $p \in P$ to $p + \langle u, p \rangle v + \langle v, p
\rangle u + a \langle u, p \rangle u$, where $a \in R$, $u, v \in P$
are fixed elements with $\langle u, v \rangle = 0$, and either $u$ or
$v$ is unimodular. Define $\sigma X$ as the map which takes $p \in P$
to either $p + \langle u, p \rangle vX + \langle vX, p \rangle u + aX
\langle u, p \rangle u$ or $p + \langle uX, p \rangle v + \langle v, p
\rangle uX + a \langle uX, p \rangle uX$. This choice depends on
whether $u$ is unimodular or $v$ is unimodular. Note that $uX$
represents $u$ {\it times} $X$, $vX$ represents $v$ {\it times} $X$,
and $a X$ represents $a$ {\it times} $X$. Also, note that $uX, vX \in
P[X]$ and $\alpha X \in R[X]$. We set $\beta(X)$ to be the product of
elements of the form $\sigma X$, whenever $\sigma$ appears in the
expression of $\alpha$. Then $\beta(1)=\alpha$ and $\beta(0)=Id$.
\hfill{$\square$}

\begin{thm} \label{eql-rel} Let $R$ be a commutative ring with $R=2R$,
  and let $I$ be an ideal of $R$. Let $(P, \langle , \rangle)$ be a
  symplectic $R$-module with $P$ a finitely generated projective module
  of rank $2n$, $n \ge 1$. Also assume that over the local ring
  $R_\gm$ for any maximal ideal $\gm$ of $R$, the alternating form
  $\langle , \rangle$ corresponds to the alternating matrix
  $\varphi_\gm$ (w.r.t. some basis), where $\varphi_\gm \equiv \psi_n ~ ({\rm
    mod}~I)$. Then $\Trans_{\Sp}(Q,IQ,\langle , \rangle) =
  \ETrans_{\Sp}(Q,IQ,\langle , \rangle)$.
\end{thm}

Proof: We have $\ETrans_{\Sp}(Q,IQ,\langle,\rangle) \subseteq
\Trans_{\Sp}(Q,IQ,\langle,\rangle)$. We need to show other inclusion.
Let us choose $\alpha$ from $\Trans_{\Sp}(Q,IQ,\langle,\rangle)$.  By
Lemma \ref{symp-trans-hom-to-identity} there exists $\alpha(X)$ in
$\Trans_{\Sp}(Q[X],IQ[X],\langle,\rangle)$ such that $\alpha(1) =
\alpha$ and $\alpha(0) = Id$. Note that
$\Trans_{\Sp}(Q_{\gm}[X],IQ_{\gm}[X],\langle,\rangle_{\psi_1 \perp
  \varphi_{\gm}}) = \ETrans_{\Sp}(Q_{\gm}[X],IQ_{\gm}[X],\langle,
\rangle_{\psi_1 \perp \varphi_{\gm}})$, for each maximal ideal $\gm$
of $R$ (follows from Remark \ref{localcase}). Hence $\alpha(X)_\gm$
belongs to
$\ETrans_{\Sp}(Q_{\gm}[X],IQ_{\gm}[X],\langle,\rangle_{\psi_1 \perp
  \varphi \otimes R_{\gm}[X]})$, for each maximal ideal $\gm$ of
$R$. Therefore, $\alpha(X) \in
\ETrans_{\Sp}(Q[X],IQ[X],\langle,\rangle)$ (see Lemma
\ref{LG-ETrans-rel}). Substituting $X=1$ we get the result.
\hfill{$\square$}

\section{\large{Equality of orbits}}

In this section we establish main result of this article regarding
equality of orbits.

\begin{thm} \label{mainTrans-rel} Let $R$ be a commutative ring with
  $R=2R$, and let $I$ be an ideal of $R$. Let $(P,\langle,\rangle)$ be
  a symplectic $R$-module with $P$ finitely generated projective
  module of rank $2n$, $n \ge 1$. Let $Q=R^2\oplus P$ and $v=(a,b,p)
  \in \Um(Q,IQ)$. We also assume for each maximal ideal $\gm$ of $R$,
  the alternating form $\langle, \rangle$ over the local ring $R_\gm$
  corresponds to an alternating matrix $\varphi_\gm$ (w.r.t. some
  basis), such that $\varphi_\gm = (1 \perp \varepsilon(\gm))^t ~
  \psi_n ~ (1 \perp \varepsilon(\gm))$, for some $\varepsilon(\gm) \in
  E_{2n-1}(R_\gm, I_\gm)$. Then
\begin{eqnarray*}
(a,b,p) ~ \ETrans(Q,IQ) & = & (a,b,p) ~ \ETrans_{\Sp}(Q,IQ, \langle,
  \rangle).
\end{eqnarray*}
\end{thm}

Proof: Let $\alpha \in \ETrans(Q,IQ)$. Let us choose $\alpha(X)$ from
$\ETrans(Q[X],IQ[X])$ such that $\alpha(1) = \alpha$ and $\alpha(0) =
Id$ (see Lemma \ref{lin-trans-hom-to-identity}). Let us define $V(X) =
(a,b,p) \alpha(X)$.

Let $\gm$ be a maximal ideal of $R$. Over $R_\gm$, we have
$\varphi_\gm = (1 \perp \varepsilon(\gm))^t ~ \psi_{n} ~ (1 \perp
\varepsilon(\gm))$, where $\varepsilon(\gm) \in \E_{2n}(R_\gm,I_\gm)$.
Let us define $W(X) = V(X) ~ (1 \perp \varepsilon(\gm))^{t}$. We have
\begin{eqnarray*}
W(X) & \in & W(0) ~ \E_{2n+2}(R_\gm[X],I_\gm[X]) \\
& = &  W(0) ~ \ESp_{2n+2}(R_\gm[X],I_\gm[X]) \\
& = &  W(0) ~ \ETrans_{\Sp}(Q_\gm[X], IQ_\gm[X]),\langle, \rangle_{\psi_{n+1}}),
\end{eqnarray*}
and hence $V(X) \in V(0) ~ \ETrans_{\Sp}(Q_\gm[X], IQ_\gm[X], \langle,
\rangle_{\psi_1 \perp \varphi_\gm})$. The first equality above follows
from (\cite{cr}, Theorem 5.6) when $n > 2$. For $n = 2$ the result has
been proved in the Appendix (see Theorem \ref{thm:relative}). This is
true for all maximal ideal $\gm$ of $R$, hence by Theorem
\ref{LG-ETrans-action-rel}, $V(X) \in V(0) ~ \ETrans_{\Sp}(Q[X],
IQ[X], \langle, \rangle)$. Substituting $X=1$ we get $ (a,b,p) ~\alpha
\in (a,b,p) ~ \ETrans_{\Sp}(Q, IQ, \langle, \rangle)$.

Now we show the other inclusion. Let $\beta \in \ETrans_{\Sp}(Q, IQ,
\langle, \rangle)$.  Let us choose $\beta(X)$ from $\ETrans_{\Sp}(Q[X],
IQ[X], \langle, \rangle)$ such that $\beta(1) = \beta$ and $\beta(0) =
Id$ (see Lemma \ref{symp-trans-hom-to-identity}). We define $V(X) =
(a,b,p) \beta(X)$.

Let $R_\gm$ be the local ring at the maximal ideal $\gm$. We
define $W(X) = V(X) ~ (1 \perp \varepsilon(\gm))^{t}$. We have
\begin{eqnarray*}
W(X) & \in & W(0) ~ ((1 \perp \varepsilon(\gm))^t)^{-1} ~ \ETrans_{\Sp}(Q_\gm[X], 
IQ_\gm[X], \langle, \rangle_{\psi_1 \perp \varphi_\gm}) ~ (1 \perp \varepsilon(\gm))^{t}\\
& = &  W(0) ~ \ETrans_{\Sp}(Q_\gm[X], IQ_\gm[X], \langle, \rangle_{\psi_{n+1}}) \\
& = &  W(0) ~ \ESp_{2n+2}(R_\gm[X],I_\gm[X]) \\
& = &  W(0) ~ \E_{2n+2}(R_\gm[X], I_\gm[X]),
\end{eqnarray*}
and hence $V(X) \in V(0) ~ \ETrans(Q_\gm[X], IQ_\gm[X])$. The last
equality above follows from (\cite{cr}, Theorem 5.6) when $n > 2$. For
$n = 2$ the result has been proved in the Appendix (see Theorem
\ref{thm:relative}). This is true for all maximal ideal $\gm$ of $R$,
hence by Theorem \ref{linear,LGaction-rel}, $V(X) \in V(0)
\ETrans(Q[X],IQ[X])$. Substituting $X=1$, we get
$(a,b,p) ~\beta \in (a,b,p) ~\ETrans(Q,IQ)$.
\hfill{$\square$}

\section{\large{Decrease in injective stability for $\Sp(Q,\langle,
  \rangle)/\ETrans_{\Sp}(Q,\langle, \rangle)$}}

Final goal of this section is to give an improvement for Basu-Rao
(\cite{BR}, Theorem 2) estimate in the module case over finitely
generated rings. For this purpose we state and prove a few preliminary
results.  {\it While dealing with the results in the relative case
  w.r.t.  an ideal $I$ of a ring $R$, we will always assume that over
  the local ring $R_\gm$, where $\gm$ is a maximal ideal, the
  alternating form $\langle, \rangle$ corresponds to the alternating
  matrix $\varphi_\gm$ (w.r.t. some basis), where $\varphi_\gm \equiv
  \psi_n ~({\rm mod}~I)$}.

\begin{thm} \label{bass-action} {\rm (\cite{bass1}, Theorem 3.4, Page
    183)} Let $R$ be a commutative ring of dim $d$. Let $I$ be an
  ideal of $R$ and $P$ be a projective module of rank $\ge d+1$.  Let
  $\tilde{Q}=R \oplus P$. Let $v_1, v_2 \in \Um(\tilde{Q})$ with $v_1
  \equiv v_2 ~({\rm mod} ~ I\tilde{Q})$. Then there exists $\beta
  \in \ETrans(\tilde{Q}, I\tilde{Q})$ such that $v_1 \beta = v_2$.
  \hfill{$\square$}
\end{thm}

\begin{thm} \label{Zalgebra} {\rm (\cite{SV}, Theorem 17.2)} Let $R$ be a 
finitely generated $\mathbb{Z}$-algebra and $kR=R$ for some natural number 
$k \ge 2$. Then $\sr(R) \le max(2, dim(R))$.
\end{thm}

\begin{thm} \label{Falgebra} {\rm (\cite{SV}, Corollary 17.3)} Let $R$ be a 
finitely generated algebra over a field which is algebraic over a finite 
field. Then $\sr(R) \le max(2, dim(R))$.
\end{thm}

Using the above three results one can prove the following theorem.

\begin{thm} \label{KumarMurthyRoy} {\rm (\cite{kmr}, Theorem 2.4)} Let $R$ be 
a finitely generated ring over $\mathbb{Z}$, or over a field $K$ which is either 
finite or algebraic closure of a finite field. Let the dimension of $R$ be $d$ 
($d \ge 2$) and $I$ be an ideal of $R$. Let $P$ be a projective module of 
rank $\ge d$ and $\tilde{Q}=R \oplus P$. Let $(a,x) \in R \oplus P$ be a unimodular 
element such that $(a, x) \equiv (1,0) ~(mod~I\tilde{Q})$. Then there exists 
$\varepsilon \in \ETrans(\tilde{Q}, I\tilde{Q})$ such that $(a, x) ~ \varepsilon 
= (1,0)$.
\hfill{$\square$}
\end{thm}

\begin{lem} \label{bass-action-Sp} Let $R$ be 
a finitely generated ring over $\mathbb{Z}$, or over a field $K$ which is either 
finite or algebraic closure of a finite field. Let the dimension of $R$ be $d$, 
and $I$ be an ideal of $R$. Let us assume $R=2R$.
  Let $(P, \langle, \rangle)$ be a symplectic $R$-module with $P$
  finitely generated projective module of even rank $\ge max \{2, d-1
  \}$, and let $Q=R^2 \oplus P$.  Let $v_1, v_2 \in \Um(Q)$ with $v_2
  = (1,0)$ and $v_1 \equiv v_2 ~({\rm mod} ~ IQ)$.  Then there exists
  $\beta \in \ETrans_{\Sp}(Q, IQ)$ such that $v_1 \beta = v_2$.
\end{lem}

Proof: Follows from Theorem \ref{KumarMurthyRoy} and Theorem
\ref{mainTrans-rel}.  \hfill{$\square$}

\begin{thm} \label{app2-trans} Let $R$ be 
a finitely generated ring over $\mathbb{Z}$, or over a field $K$ which is either 
finite or algebraic closure of a finite field. Let the dimension of $R$ be $d$. Let 
us assume $R=2R$. Let $(P, \langle, \rangle)$ be a
  symplectic $R$-module with $P$ finitely generated projective module
  of even rank $\ge d-3 $.  Let $Q=(R^2 \oplus P)$, and let
  $\widehat{Q}=(R^2 \oplus Q)$.  Let $\sigma \in \Sp(Q, \langle,
  \rangle)$ with $(I_2 \perp \sigma) \in \ETrans_{\Sp}(\widehat{Q},
  \langle, \rangle )$.  Then $\sigma$ is (stably elementary
  symplectic) homotopic to the identity. In fact, $\sigma = \rho(1)$,
  with $\rho(0) = Id.$, for some $\rho(X) \in \Sp(Q[X],\langle,
  \rangle)$ such that $(I_2 \perp \rho(X)) \in
  \ETrans_{\Sp}(\widehat{Q}[X], \langle, \rangle)$.
\end{thm}

Proof: Let us choose $\alpha(X)$ from $\ETrans_{\Sp}(\widehat{Q}[X],
\langle, \rangle)$, such that $\alpha(1) = (I_2 \perp \sigma)$, and
$\alpha(0)=Id.$ (see Lemma \ref{symp-trans-hom-to-identity}). Let
$v(X) = e_1 \alpha(X)$. Note that $v(X)$ belongs to
$\Um(\widehat{Q}[X],(X^2 - X) \widehat{Q}[X])$.  Also, $e_1 =(1,0,0)$
belongs to $\Um(\widehat{Q}[X],(X^2 - X) \widehat{Q}[X])$. Therefore,
by Lemma \ref{bass-action-Sp}, we have $\beta(X) \in
\ETrans_{\Sp}(\widehat{Q}[X], (X^2 - X)\widehat{Q}[X], \langle,
\rangle)$, such that $v(X) \beta(X) =(1,0,0)$, i.e, $e_1 \alpha(X)
\beta(X) = e_1$. Let $\delta(X)$ denote the product $\alpha(X)
\beta(X)$. Note that $\delta(X) \in \ETrans_{\Sp}(\widehat{Q}[X],
\langle, \rangle)$ and $e_1 \delta(X) = e_1$. Therefore, $\delta(X)$ will be of the form
\begin{eqnarray*}
\begin{pmatrix} 1 & 0 & 0 \\ * & 1 & *
  \\  * & 0 & \gamma(X) \end{pmatrix},  
\end{eqnarray*}
where $\gamma(X) \in \Sp(Q[X], \langle, \rangle)$. Note that
$\delta(1) = (I_2 \perp \sigma)$. Therefore $\gamma(1)=\rho$,
$\gamma(0)=Id.$, and hence $\rho$ is (symplectic) homotopic to
identity.
\hfill{$\square$}

\begin{thm} \label{vorst-trans} {\rm (\cite{BR}, Theorem 3.13)}
Let $(R,\gm)$ be a regular local ring. Assume that $R$ contains a
field, or characteristic of $R / \gm$ is not in $\gm^2$. Then
\begin{eqnarray*}
\Sp(R[X]^{2n+2}, \langle, \rangle_{\varphi_\gm}) & = & \ETrans_{\Sp}(R[X]^{2n+2},
\langle, \rangle_{\varphi_\gm}),
\end{eqnarray*}
for $n \ge 1$, where $\varphi_\gm$ is the associated matrix of the
alternating bilinear form $\langle, \rangle$ (w.r.t. some basis).
\hfill{$\square$}
\end{thm}

The next corollary improves Basu-Rao estimate in the module
case over finitely generated rings.

\begin{cor} \label{app3-trans} Let $R$ be a finitely generated
  non-singular algebra of dimension $d$ over $K$, where $K$ is either
  $\mathbb{Z}$, or a finite field or the algebraic closure of a finite field. Let us
  assume $R=2R$. Let $(P, \langle, \rangle)$ be a symplectic
  $R$-module with $P$ finitely generated projective module of even
  rank $\ge max \{2, d-3 \}$.  Let $Q=(R^2 \oplus P)$, and
  $\widehat{Q}=(R^2 \oplus Q)$. If $\sigma \in \Sp(Q, \langle,
  \rangle)$ with $(I_2 \perp \sigma) \in \ETrans_{\Sp} (\widehat{Q},
  \langle, \rangle )$, then $\sigma$ belongs to $\ETrans_{\Sp} (Q,
  \langle, \rangle )$.
\end{cor}

Proof: From the proof of Theorem \ref{app2-trans} it follows that
$\sigma = \rho(1)$ for some $\rho(X) \in \Sp(Q[X],\langle,\rangle)$,
with $\rho(0) = Id$. Using Theorem \ref{vorst-trans} we get that
$\Sp(R_\gm[X]^{2n+2}, \langle,
\rangle_{\varphi_\gm})\\ =~\ETrans_{\Sp}(R_\gm[X]^{2n+2}, \langle,
\rangle_{\varphi_\gm})$, for all maximal ideals $\gm$ of $R$. This
implies $\rho(X) \in \ETrans_{\Sp}(R_\gm[X]^{2n+2}, \langle,
\rangle_{\varphi_\gm})$, for all maximal ideals $\gm$ in
$R$. Therefore, $\rho(X) \in \ETrans_{\Sp}(R[X], \langle, \rangle)$
(see Lemma \ref{LG-ETrans-rel}). Hence $\sigma = \rho(1)$ belongs to
$\ETrans_{\Sp}(R, \langle, \rangle)$.  \hfill{$\square$}

\bigskip

\begin{center}
 {\bf \Large Appendix} \\[5mm]
{\bf Alpesh Dhorajia}\\
\small{Birla Institute of Technology and Science,
K. K. Birla Goa Campus,
Goa 403 726, India}
\end{center}

The aim of this Appendix is to work out via commutator laws the size
$4$ cases in the study above. This entails proving Lemma
\ref{vanderk1} and Lemma 3.1 above in the symplectic size $4$ case.

In order to prove Lemma \ref{vanderk1} in the symplectic case we need
to establish the following inclusion first.

\begin{lem} \label{inclusion} 
{\rm Let $n \ge 2$, $R$ be a commutative ring with
  $R=2R$, and let $I$ be an ideal of $R$.  Then $\ESp_{2n}(R,I)
  \subseteq \ESp_{2n}^1(R,I)$.}
\end{lem}

Proof: It suffices to show $\ESp_{2n}^1(R,I)$ contains the set $S_{ij}
= \{ se_{ij}(a) se_{ji}(x) se_{ij}(-a): a \in R, x \in I \}$, for all
$i,j$, with $i \ne j$. Note that $S_{ij}=S_{\sigma(j) \sigma(i)}$ when
$i\ne \sigma(j)$. First we state the following identities
\begin{align}\label{commprop}
[gh,k] &~ = ~ \big({}^g[h,k]\big)[g,k],\\
[g,hk] &~ = ~ [g,h]\big({}^h[g,k]\big),\\
{}^g[h,k] &~ = ~ [{}^gh,{}^gk],
\end{align}
where $^gh$ denotes $ghg^{-1}$ and $[g,h]=ghg^{-1}h^{-1}$. These
identities will be used throughout the proof without mentioning it
always.

Also, we need the following relations which hold for all integers $i,j$
with $i \neq j, \sigma(j)$, and for all $a, b
\in R$.
\begin{align}\label{commprop2}
[se_{i \sigma(i)}(a), se_{\sigma(i) j}(b)] ~=~ se_{ij}(ab)
se_{\sigma(j) j}((-1)^{i+j}ab^2), \\ [se_{ik}(a), se_{kj}(b)] ~=~
se_{ij}(ab), {\rm if} ~ k \ne \sigma(i), \sigma(j), \\ [se_{ik}(a),
  se_{k \sigma(i)}(b)] ~=~ se_{i \sigma(i)}(2ab), {\rm if} ~ k \ne i,
\sigma(i).
\end{align}
The following relation holds for all $i, j$ with $1 \le i\ne j \le
2n$, and for all $a, b \in R$.
\begin{align}
[se_{ij}(a), se_{kl}(b)] ~=~ Id., {\rm if} ~ i \ne l, \sigma(k) ~ {\rm
  and} ~ j \ne k, \sigma(l).
\end{align}

Note that $se_{1j}(x), se_{i1}(x) \in \ESp_{2n}^{1}(R,I)$, for $2 \le
i,j \le 2n$ and $x \in I$. For $3 \le i,j \le 2n$ with $i \ne
\sigma(j)$, and $x \in I$, we have $se_{ij}(x) = [se_{i1}(x),
  se_{1j}(1)] \in \ESp_{2n}^{1}(R,I)$ and $se_{i \sigma(i)}(x) =
      [se_{i1}(x/2), se_{1 \sigma(i)}(1)] \in \ESp_{2n}^{1}(R,I)$. In
      the following computation we will express the generators of
      $S_{ij}$ in terms of $se_{ij}(x)$ or $se_{1j}(a)$, where $x \in
      I, a \in R $. Also, note that we will use $*$ to represent
      elements of the ideal $I$ and $\circledast$ to represent
      elements of the ring $R$, which does not belong to the ideal
      $I$.

{\it Case (a):} It is clear that $S_{1j} \subseteq \ESp_{2n}^{1}(R,I)$,
for $2 \le j \le 2n$.

{\it Case (b):} Here we show $S_{i1} \subseteq \ESp_{2n}^{1}(R,I)$,
for $3 \le i \le 2n$. Let $a \in R$ and $x \in I$. We have
\begin{eqnarray*}
&& ^{se_{ij}(a)} se_{ji}(x) \\ &=& ^{se_{ij}(a)} ([se_{j
      \sigma(j)}(1), se_{\sigma(j) i}(*)] se_{\sigma(i) i}(*))\\ &=&
  [se_{j \sigma(i)}(\circledast) se_{i \sigma(i)}(\circledast) se_{j
      \sigma(j)}(1), se_{\sigma(j) j}(*) se_{\sigma(j) i}(*)]
  ^{se_{ij}(a)} se_{\sigma(i) i}(*)\\ &=& [se_{i
      \sigma(i)}(\circledast) se_{j \sigma(i)}(\circledast) se_{j
      \sigma(j)}(1), se_{\sigma(j) j}(*) se_{\sigma(j) i}(*)]
  ^{se_{ij}(a)} se_{\sigma(i) i}(*)\\ &=& \big( ^{se_{i
      \sigma(i)}(\circledast)} [se_{j \sigma(i)}(\circledast) se_{j
      \sigma(j)}(1), se_{\sigma(j) j}(*) se_{\sigma(j) i}(*)]
  \big)\\ && [se_{i \sigma(i)}(\circledast), se_{\sigma(j) j}(*)
    se_{\sigma(j) i}(*)] ^{se_{ij}(a)} se_{\sigma(i) i}(*)\\ &=& ABC
  ~ {\rm (say)},
\end{eqnarray*}
where
\begin{eqnarray*}
A &=&  ^{se_{i
      \sigma(i)}(\circledast)} [se_{j \sigma(i)}(\circledast) se_{j
      \sigma(j)}(1), se_{\sigma(j) j}(*) se_{\sigma(j) i}(*)] ,\\
B &=& [se_{i \sigma(i)}(\circledast), se_{\sigma(j) j}(*) se_{\sigma(j)
      i}(*)],\\ 
C &=& ^{se_{ij}(a)} se_{\sigma(i) i}(*).
\end{eqnarray*}

Now
\begin{eqnarray*}
A &=& ^{se_{i \sigma(i)}(\circledast)} [se_{j
    \sigma(i)}(\circledast) se_{j \sigma(j)}(1), se_{\sigma(j) j}(*)
  se_{\sigma(j) i}(*)], \\ &=& [se_{j \sigma(i)}(\circledast)
  se_{j \sigma(j)}(1), se_{\sigma(j) j}(*) se_{ij}(*) se_{\sigma(j)
    j}(*) se_{\sigma(j) i}(*)] \\ &=& [se_{j \sigma(j)}(1) se_{j
    \sigma(i)}(\circledast) , se_{\sigma(j) j}(*) se_{ij}(*)
  se_{\sigma(j) i}(*)],\\ B &=& [se_{i \sigma(i)}(\circledast),
  se_{\sigma(j) j}(*) se_{\sigma(j) i}(*)] \\ &=& ^{se_{i
    \sigma(i)}(\circledast)} se_{\sigma(j) j}(*) ^{se_{i
    \sigma(i)}(\circledast)} se_{\sigma(j) i}(*) se_{\sigma(j)
  i}(*)^{-1} se_{\sigma(j) j}(*)^{-1} \\ &=& se_{\sigma(j) j}(*)
se_{ij}(*) se_{\sigma(j) j}(*) se_{\sigma(j) i}(*) se_{\sigma(j) i}(*)
se_{\sigma(j) j}(*), \\ C &=& ^{se_{ij}(a)} se_{\sigma(i) i}(*))
\\ &=& se_{\sigma(i) j}(*) se_{\sigma(j) j}(*) se_{\sigma(i) i}(*)).
\end{eqnarray*}

It is clear that $B, C \in \ESp_{2n}^{1}(R,I)$. In this case we have
$j=1$, and hence $A \in \ESp_{2n}^{1}(R,I)$.

{\it Case (c):} Here we show $S_{ij} \subseteq \ESp_{2n}^{1}(R,I)$,
for $3 \le i, j \le 2n$, and $j=\sigma(i)$. Note that 
\begin{eqnarray*}
&& ^{se_{i \sigma(i)}(a)} se_{\sigma(i) i}(x) \\ &=& ^{se_{i
      \sigma(i)}(a)} [se_{\sigma(i) 1}(x/2), se_{1i}(1)]\\ &=&
  [se_{i1}(*) se_{21}(*) se_{\sigma(i) 1}(*), se_{i2}(\circledast)
    se_{12}(\circledast) se_{1i}(1)]\\ &=& [se_{i1}(*) se_{21}(*)
    se_{\sigma(i) 1}(*), se_{1 \sigma(i)} (\circledast)
    se_{12}(\circledast) se_{1i}(1)] \\ &\in& \ESp_{2n}^{1}(R,I).
\end{eqnarray*}

{\it Case (d):} Here we show $S_{ij} \subseteq \ESp_{2n}^{1}(R,I)$, for
$3 \le i, j \le 2n$, and $j \ne \sigma(i)$. From the computaion in
{\it Case (b)} it is clear that $^{se_{ij}(a)} se_{ji}(x)$ can be
expressed as product of $ABC$, where $B, C \in
\ESp_{2n}^{1}(R,I)$. Note that
\begin{eqnarray*}
A &=& [se_{j \sigma(j)}(1) se_{j \sigma(i)}(\circledast) ,
  se_{\sigma(j) j}(*) se_{ij}(*) se_{\sigma(j) i}(*)] \\ &=& ^{se_{j
    \sigma(j)}(1) se_{j \sigma(i)}(\circledast)} se_{\sigma(j) j}(*) ~
^{se_{j \sigma(j)}(1) se_{j \sigma(i)}(\circledast)} se_{ij}(*) ~
^{se_{j \sigma(j)}(1) se_{j \sigma(i)}(\circledast)} se_{\sigma(j)
  i}(*) \\ && se_{\sigma(j) i}(*)^{-1} se_{ij}(*)^{-1} se_{\sigma(j)
  j}(*)^{-1} \\ &=& PQR ~ se_{\sigma(j) i}(*) se_{ij}(*) se_{\sigma(j)
  j}(*),
\end{eqnarray*}
where
\begin{eqnarray*}
P &=& ^{se_{j \sigma(j)}(1) se_{j \sigma(i)}(\circledast)}
se_{\sigma(j) j}(*), \\ Q &=& ^{se_{j \sigma(j)}(1) se_{j
    \sigma(i)}(\circledast)} se_{ij}(*), \\ R &=& ^{se_{j
    \sigma(j)}(1) se_{j \sigma(i)}(\circledast)} se_{\sigma(j) i}(*).
\end{eqnarray*}

Now
\begin{eqnarray*}
P &=& ^{se_{j \sigma(j)}(1) se_{j \sigma(i)}(\circledast)}
se_{\sigma(j) j}(*) \\ &=& ^{se_{j \sigma(j)}(1)} \big( se_{\sigma(j)
  \sigma(i)}(*) se_{i \sigma(i)}(*) se_{\sigma(j) j}(*) \big) \\ &=&
se_{j \sigma(i)}(*) se_{i \sigma(i)}(*) se_{\sigma(j) \sigma(i)}(*)
se_{i \sigma(i)}(*) \big( ^{se_{j \sigma(j)}(1)} se_{\sigma(j) j}(*) \big),
\end{eqnarray*}
\begin{eqnarray*}
Q &=& ^{se_{j \sigma(j)}(1) se_{j \sigma(i)}(\circledast)} se_{ij}(*)
\\ &=& ^{se_{j \sigma(j)}(1)} se_{i \sigma(i)}(*) \\ &=& se_{i
  \sigma(i)}(*), \\ R &=& ^{se_{j \sigma(j)}(1) se_{j
    \sigma(i)}(\circledast)} se_{\sigma(j) i}(*) \\ &=& ^{se_{j
    \sigma(i)}(\circledast) se_{j \sigma(j)}(1)} se_{\sigma(j) i}(*)
\\ &=& ^{se_{j \sigma(i)}(\circledast)} \big( se_{ji}(*) se_{\sigma(i)
  i}(*) se_{\sigma(j) i}(*) \big) \\ &=& ^{se_{j
    \sigma(i)}(\circledast)} \big( se_{ji}(*) se_{\sigma(i) i}(*)
se_{\sigma(i) j}(*) \big) \\ &=& se_{j \sigma(j)}(*) se_{ji}(*)
se_{\sigma(i) \sigma(j)}(*) se_{j \sigma(j)}(*) se_{\sigma(i) i}(*)
\\ && ^{se_{j \sigma(i)}(\circledast)} se_{\sigma(i) j}(*) \\
&=& se_{j \sigma(j)}(*) se_{ji}(*) se_{ji}(*) se_{j \sigma(j)}(*)
  se_{\sigma(i) i}(*) \\ && ^{se_{j \sigma(i)}(\circledast)}
  [se_{\sigma(i) 1}(*), se_{1j}(1)] \\
&=& se_{j \sigma(j)}(*)
  se_{ji}(*) se_{j \sigma(j)}(*) se_{\sigma(i) i}(*)\\ && [se_{j1}(*)
    se_{\sigma(i) 1}(*), se_{1 \sigma(i)}(\circledast) se_{1j}(1)]. \\
\end{eqnarray*}
and hence $A \in \ESp_{2n}^{1}(R,I)$.

{\it Case (e):} Finally we show $S_{21} \subseteq \ESp_{2n}^{1}(R,I)$. Note that
\begin{eqnarray*}
&& ^{se_{21}(a)} se_{12}(x) \\ &=& ^{se_{21}(a)} [se_{1k}(1),
    se_{k2}(*)], ~{\rm where}~k \ne 1,2 \\ &=& [^{se_{21}(a)}
    se_{1k}(1), se_{2 \sigma(k)}(*) se_{k \sigma(k)}(*)
    se_{k2}(*)]\\ &=& [\alpha , se_{2 \sigma(k)}(*) se_{k
      \sigma(k)}(*) se_{k2}(*)], ~{\rm where}~ \alpha = ^{se_{21}(a)}
  se_{1k}(1) \\ &=& ^{\alpha} se_{2 \sigma(k)}(*) ^{\alpha} se_{k
    \sigma(k)}(*) ^{\alpha} se_{k2}(*) se_{k2}(*)^{-1} se_{k
    \sigma(k)}(*)^{-1} se_{2 \sigma(k)}(*)^{-1}. \\
\end{eqnarray*} 

The $\alpha$ above can be expressed in the following three different
ways, namely
\begin{eqnarray*}
\alpha &=& se_{1k}(1) se_{2k}(\circledast) se_{\sigma(k)
  k}(\circledast) \\ &=& se_{2k}(\circledast) se_{\sigma(k) k}(\circledast)
se_{1k}(1) \\ &=& se_{\sigma(k) k}(\circledast) se_{2k}(\circledast) se_{1k}(1),
\end{eqnarray*} 
and we use them as we find it convenient. Now 
\begin{eqnarray*}
^{\alpha} se_{2 \sigma(k)}(*) &=& ^{\alpha} se_{k1}(*) \\ &=&
  ^{se_{1k}(1) se_{2k}(\circledast) se_{\sigma(k) k}(\circledast)}
  se_{k1}(*) \\ &=& ^{se_{1k}(1) se_{2k}(\circledast)} \big(
  se_{\sigma(k) 1}(*) se_{21}(*) se_{k1}(*) \big) \\ &=& ^{se_{1k}(1)}
  \big( se_{\sigma(k) 1}(*) se_{21}(*) \big) \\ &=& se_{\sigma(k)
    k}(*) se_{\sigma(k) 1}(*) se_{2k}(*) se_{\sigma(k) k}(*)
  se_{21}(*) \\ &=& se_{\sigma(k) k}(*) se_{2k}(*) se_{21}(*),
\end{eqnarray*} 
\begin{eqnarray*}  
^{\alpha} se_{k2}(*) &=& ^{se_{2k}(\circledast) se_{\sigma(k)
      k}(\circledast) se_{1k}(1)} se_{k2}(*) \\ &=&
  ^{se_{2k}(\circledast) se_{\sigma(k) k}(\circledast)} \big(
  se_{12}(*) se_{k2}(*) \big) \\ &=& ^{se_{2k}(\circledast)} \big(
  se_{12}(*) se_{\sigma(k) 2}(*) se_{12}(*) se_{k2}(*) \big) \\ &=&
  ^{se_{2k}(\circledast)} \big( se_{12}(*) se_{\sigma(k) 2}(*)
  se_{k2}(*) \big) \\ &=& se_{1k}(*) se_{\sigma(k) k}(*) se_{12}(*)
  se_{\sigma(k) k}(*) se_{\sigma(k) 2}(*) \big( ^{se_{2k}(\circledast)}
  se_{k2}(*) \big), \\ 
^{\alpha} se_{k \sigma(k)}(*) &=&
  ^{se_{\sigma(k) k}(\circledast) se_{2k}(\circledast) se_{1k}(1)}
  se_{k \sigma(k)}(*)\\ &=& ^{se_{\sigma(k) k}(\circledast)
    se_{2k}(\circledast)} \big( se_{k2}(*) se_{12}(*) se_{k
    \sigma(k)}(*) \big) \\ &=& XYZ ~{\rm(say)},
\end{eqnarray*} 
where
\begin{eqnarray*}
X &=& ^{se_{\sigma(k) k}(\circledast) se_{2k}(\circledast)} se_{k2}(*)
\\ &=& ^{se_{2k}(\circledast) se_{\sigma(k) k}(\circledast)}
se_{k2}(*) \\ &=& ^{se_{2k}(\circledast)} \big( se_{\sigma(k) 2}(*)
se_{12}(*) se_{k2}(*) \big) \\ &=& se_{\sigma(k) k}(*) se_{\sigma(k)
  2}(*) se_{1k}(*) se_{\sigma(k) k}(*) se_{12}(*) ~
\big( ^{se_{2k}(\circledast)} se_{k2}(*) \big), \\ Y &=& ^{se_{\sigma(k)
    k}(\circledast) se_{2k}(\circledast)} se_{12}(*) \\ &=&
^{se_{\sigma(k) k}(\circledast)} \big( se_{1k}(*) se_{\sigma(k) k}(*)
se_{12}(*) \big) \\ &=& se_{1k}(*) se_{\sigma(k) k}(*) se_{12}(*), \\ Z
&=& ^{se_{\sigma(k) k}(\circledast) se_{2k}(\circledast)} se_{k
  \sigma(k)}(*) \\ &=& ^{se_{\sigma(k) k}(\circledast)} \big(
se_{k1}(*) se_{21}(*) se_{k \sigma(k)}(*) \big) \\ &=& se_{\sigma(k)
  1}(*) se_{21}(*) se_{k1}(*) se_{21}(*) ~ \big( ^{se_{\sigma(k)
    k}(\circledast)} se_{k \sigma(k)}(*) \big).
\end{eqnarray*} 
\hfill{$\square$}

\medskip
Proof of Lemma \ref{vanderk1}:
Let $f : \E^1(n,R,I) \lra \E^1(n,R/I,0)$. Note that 
\begin{eqnarray*}
\ker(f) &\subseteq& \E^1(n,R,I) \cap \G(n,R,I).
\end{eqnarray*}

Let $M=\prod ge_{j1}(x_j) ge_{1i}(a_i) \in \E^1(n,R,I) \cap
\G(n,R,I)$. Note that $M \in \G(n,R,I)$ implies $\overline{M} = \prod
ge_{j1}(0)ge_{1i}(\bar{a_i}) = \prod ge_{1i}(\bar{a_i}) = Id$. Here
{\it bar} means reduction modulo the ideal $I$. Therefore, $M \in
\ker(f)$ and hence $\ker(f) = \E^1(n,R,I) \cap \G(n,R,I)$.

Now we shall prove that $\ker(f) = \E(n,R,I)$. Let $E = \prod_{k=1}^r
ge_{j_k 1}(x_k) ge_{1i_k}(a_k) \in \ker(f)$.  Note that $E$ can be
written as $ge_{j_1 1}(x_1) \prod_{k=2}^r \gamma_k ge_{j_k 1}(x_k)
\gamma_k^{-1}$, where $\gamma_l$ is equal to
$\prod_{k=1}^{l-1}~ge_{1i_k}(a_k) \in \E(n,R)$, and hence $\ker(f)
\subseteq \E(n,R,I)$.  The reverse inclusion follows from the fact
that $\E(n,R,I) \subseteq \E^1(n,R,I)$.  \hfill{$\square$}

\medskip
Proof of Lemma \ref{ness4-E1} in the symplectic case when $n=4$: Given
that $se_{ij}(Xf(X))$ is in $\ESp_{2n}^1(R[X],I[X])$. First assume
$i=1$ and $f(X) \in R[X]$. We prove the result using induction on $r$,
the number of generators of $\varepsilon$. Let $r=1$ and $\varepsilon
= se_{pq}(a)$. Note that when $p=1, a \in R$, and when $q=1, a \in
I$. Also, note that we will use $*$ to represent elements of the
ideal $I$. Computation done in {\it Case (1)} to {\it Case (11)} in
Lemma \ref{ness4-E1} will be the same here except {\it Case (7)}. We
work out the details of {\it Case (7)} now.

{\it Case $($7$)$:} Let $(p,q)=(j,1)$ and $j \ne 2$. In this case
\begin{eqnarray*}
&& se_{j1}(a) ~ se_{1j}(Y^4X f(Y^4X)) ~ se_{j1}(-a) \\
&=& se_{j1}(a) ~ se_{\sigma(j) 2}(Y^4X f(Y^4X)) ~ se_{j1}(-a) \\
&=& ^{se_{j1}(a)} ~ \Big([ se_{\sigma(j) j}(Y^2), ~ se_{j 2}(Y^2X f(Y^4X))] ~ se_{12}((-1)^j Y^6X^2 f^2(Y^4X)) \Big) \\
&=& [se_{\sigma(j) 1}(* Y^2) se_{21}(* Y^2) se_{\sigma(j) j}(Y^2), se_{j \sigma(j)}(* Y^2X f(Y^4X)) se_{j 2}(Y^2X f(Y^4X))] \\ && se_{1 \sigma(j)}(* Y^6X^2 f^2(Y^4X)) ~ se_{j \sigma(j)}(* Y^6X^2 f^2(Y^4X)) ~ se_{12}((-1)^j Y^6X^2 f^2(Y^4X)) \\
&=& se_{\sigma(j) 1}(* Y^2) ~ se_{21}(* Y^2)  \Big( se_{\sigma(j) j}(Y^2) se_{j \sigma(j)}(* Y^2X f(Y^4X)) se_{j 2}(Y^2X f(Y^4X)) \\
&& se_{\sigma(j) j}(- Y^2) \Big) ~ se_{21}(* Y^2) ~ se_{\sigma(j) 1}(* Y^2) ~ se_{j 2}(-Y^2X f(Y^4X)) \\
&& se_{j \sigma(j)}(* Y^2X f(Y^4X)) ~ se_{1 \sigma(j)}(* Y^6X^2 f^2(Y^4X)) ~ se_{j \sigma(j)}(* Y^6X^2 f^2(Y^4X)) \\
&& se_{12}((-1)^j Y^6X^2 f^2(Y^4X)) \\
&=& se_{\sigma(j) 1}(* Y^2) ~ se_{21}(* Y^2)  \Big( ^{se_{\sigma(j) j}(Y^2)} [se_{j 1}(*Y), ~ se_{1 \sigma(j)}(* YX f(Y^4X))] \\
&& ^{se_{\sigma(j) j}(Y^2)} se_{j 2}(Y^2X f(Y^4X)) \Big) ~ se_{21}(* Y^2) ~ se_{\sigma(j) 1}(* Y^2) ~ se_{j 2}(-Y^2X f(Y^4X)) \\
&& se_{j \sigma(j)}(* Y^2X f(Y^4X)) ~ se_{1 \sigma(j)}(* Y^6X^2 f^2(Y^4X)) ~ se_{j \sigma(j)}(* Y^6X^2 f^2(Y^4X)) \\
&& se_{12}((-1)^j Y^6X^2 f^2(Y^4X)) \\
&=& se_{\sigma(j) 1}(* Y^2) ~ se_{21}(* Y^2) ~ [se_{\sigma(j) 1}(* Y^3) ~ se_{21}(* Y^4) ~ se_{j 1}(*Y), \\
&& se_{\sigma(j) 2}(* Y^3X f(Y^4X)) ~ se_{12}(* Y^4X^2 f^2(Y^4X)) ~ se_{1 \sigma(j)}(* YX f(Y^4X))] \\
&& se_{\sigma(j) 2}(Y^4X f(Y^4X)) ~ se_{12}(Y^6X^2 f^2(Y^4X)) ~ se_{j 2}(Y^2X f(Y^4X)) ~ se_{21}(* Y^2) \\
&& se_{\sigma(j) 1}(* Y^2) ~ se_{j 2}(-Y^2X f(Y^4X)) ~ [se_{j 1}(*Y), ~ se_{1 \sigma(j)}(Y X f(Y^4X))] \\
&& se_{1 \sigma(j)}(* Y^6X^2 f^2(Y^4X)) ~ [se_{j 1}(*Y), ~ se_{1 \sigma(j)}(Y^5 X^2 f^2(Y^4X))] \\
&& se_{12}((-1)^j Y^6X^2 f^2(Y^4X)) 
\end{eqnarray*}

Arguing as in the proof of Lemma \ref{ness4-E1} of \S 3, we
can complete this proof.
\hfill{$\square$}

\medskip

Using equations (6) to (9) we can improve Theorem 5.6 of
\cite{cr}. For this purpose we first establish the following lemma
(\cite{cr}, Lemma 5.2) for size 4:

\begin{lem} \label{subset-imprv} 
Let $I$ be an ideal of $R$. Assume that $2R=R$.
Then ${\rm ESp}_{2n}(R, I^2) \subseteq {\rm ESp}_{2n}(I)$, for $n \ge
    2$.
\end{lem}

Proof: Let $z^* = \sum a_t b_t$ with $a_t, b_t \in I$.  Let $\beta =
se_{i j}(z^*) \in {\rm ESp}_{2n}(I^2)$ and $\alpha= se_{k l}(z) \in
{\rm ESp}_{2n}(R)$ for some $z \in R$. It suffices to show that
$\alpha \beta \alpha^{-1} \in {\rm ESp}_{2n}(I)$. If $(k,l) \ne (j,i)$
and $(k,l) \ne (\sigma(i), \sigma(j))$, then the matrix $\alpha \beta
\alpha^{-1}$ splits into a product of elementary matrices from $ {\rm
  ESp}_{2n}(I)$. Now we consider the following two cases.

{\it Case $($1$)$:} Let $(k,l) = (j,i)$
or $(k,l) = (\sigma(i), \sigma(j))$, and $j \ne \sigma(i)$.

\begin{eqnarray*}
\alpha \beta \alpha^{-1} &=& \prod se_{ji}(z) se_{ij}(a_t b_t) se_{ji}(-z) \\
&=& \prod ~~ ^{se_{ji}(z)} \big( [se_{i \sigma(i)}(a_t), ~ se_{\sigma(i) j}(b_t)] ~ se_{\sigma(j) j}((-1)^{i+j+1} a_t b_t^2) \big)\\
&=& \prod [se_{i \sigma(j)}(-a_t z) ~ se_{j \sigma(j)}((-1)^{i+j} a_t z^2) ~ se_{i \sigma(i)}(a_t), ~ se_{\sigma(i) i}(2 b_t z) ~ se_{\sigma(i) j}(b_t)] \\
&& se_{\sigma(j) i}((-1)^{i+j} a_t b_t^2 z) ~ se_{\sigma(i) i}(-a_t b_t^2 z^2) ~ se_{\sigma(j) j}((-1)^{i+j+1} a_t b_t^2).
\end{eqnarray*}

{\it Case $($2$)$:} Let $(k,l) = (j,i)$ and $j = \sigma(i)$. Let us
choose an integer $d \ne i, \sigma(i)$.
\begin{eqnarray*}
\alpha \beta \alpha^{-1} &=& \prod se_{\sigma(i) i}(z) se_{i \sigma(i)}(a_t b_t) se_{\sigma(i) i}(-z) \\
&=& \prod ~~ ^{se_{\sigma(i) i}(z)} [se_{i d}(a_t), se_{d \sigma(i)}(b_t/2)] \\
&=& \prod [se_{\sigma(i) d}(a_t z) ~ se_{\sigma(d) d}((-1)^{i+d+1} a_t^2 z) ~ se_{i d}(a_t), se_{d i}(b_t z/2) se_{d \sigma(i)}(b_t/2)].
\end{eqnarray*}
\hfill{$\square$}

\medskip
Now we have Lemma 5.1 to Lemma 5.4 and Theorem 5.5 of \cite{cr} for
Elementary symplectic group of size 4, which will prove the size $4$ case of 
(\cite{cr}, Theorem 5.6), {\it viz.}:

\begin{thm} \label{thm:relative} {\rm (\cite{cr}, Theorem 5.6)} 
Let $I$ be an ideal of $R$. Assume $2R=R$. Then the natural map ${\rm
  Um}_{2n}(R, I)/{\rm ESp}_{2n}(R, I)\longrightarrow {\rm Um}_{2n}(R,
I)/{\rm E}_{2n}(R, I)$ is bijective for $n \ge 2$.
\hfill{$\square$}
\end{thm}

\end{document}